\documentclass[reqno,11pt]{amsart}
\usepackage{amsfonts,amsmath,amssymb,fullpage,array,graphicx,rotate}

\theoremstyle{plain} \newtheorem{Thm}{Theorem}[section]
\theoremstyle{plain} 
\theoremstyle{plain} 
\theoremstyle{plain} \newtheorem{Lemma}[Thm]{Lemma}
\theoremstyle{plain} 
\theoremstyle{definition} \newtheorem{Def}[Thm]{Definition}
\theoremstyle{remark} 
\theoremstyle{definition} \newtheorem{Ex}[Thm]{Example}

\newcommand{\clearemptydoublepage}
             {\newpage{\pagestyle{empty}\cleardoublepage}}

\newcommand{\thmlist}%
          {\renewcommand{\theenumi}{\arabic{enumi}}
                  \renewcommand{\labelenumi}{{\rm (\theenumi)}}
                            }
\newcommand{\factlist}%
          {\renewcommand{\theenumi}{\arabic{enumi}}
                  \renewcommand{\labelenumi}{{\rm (\theenumi)}}
                            }
\newcommand{\factlisti}%
          {
                  
                            }
\newcommand{\condlisti}%
          {
                  \renewcommand{\labelenumi}{{\rm(\theenumii)}}
                            }
\newcommand{\condlist}%
          {\renewcommand{\theenumi}{\roman{enumi}}
                  \renewcommand{\labelenumi}{{\rm (\theenumi)}}
                            }

\def\frak{\mathfrak}
\renewcommand{\Re}{\mathop{\rm{Re}}}

\newcommand{\conv}{\mathop{\rm{conv}}}

\newcommand{\supp}{\mathop{\rm{supp}}}
\newcommand{\Ad}{\mathop{\rm{Ad}}}
\newcommand{\ad}{\mathop{\rm{ad}}}

\newcommand{\PW}{\mathop{\rm{PW}}}
\newcommand{\abs}[1]{\left|\/#1\/\right|}

\newcommand{\inner}[2]{\langle#1,#2\rangle}

\newcommand{\C}{\ensuremath{\mathbb C}}

\newcommand{\D}{\ensuremath{\mathbb D}}

\newcommand{\R}{\ensuremath{\mathbb R}}
\newcommand{\Z}{\ensuremath{\mathbb Z}}
\newcommand{\N}{\ensuremath{\mathbb N}}
\newcommand{\quat}{\ensuremath{\mathbb H}}

\newcommand{\fa}{\ensuremath{\mathfrak a}}
\newcommand{\fg}{\ensuremath{\mathfrak g}}
\newcommand{\fh}{\ensuremath{\mathfrak h}}
\newcommand{\fk}{\ensuremath{\mathfrak k}}
\newcommand{\fn}{\ensuremath{\mathfrak n}}
\newcommand{\fm}{\ensuremath{\mathfrak m}}
\newcommand{\fp}{\ensuremath{\mathfrak p}}
\newcommand{\fq}{\ensuremath{\mathfrak q}}
\renewcommand{\l}{\lambda}
\renewcommand{\a}{\alpha}
\renewcommand{\b}{\beta}
\newcommand{\e}{\varepsilon}

\newcommand{\la}{\l_\a}
\newcommand{\G}{\Gamma}

\newcommand\SL{\mathop{\rm{SL}}}

\newcommand{\frakacs}{\frak a_{\scriptscriptstyle{\C}}^*}

\newcommand{\liecomplex}[1]{{\frak #1}_{\scriptscriptstyle{\C}}}
\newcommand{\complex}[1]{{#1}_{\scriptscriptstyle{\C}}}

\newcommand{\rootstheta}{\langle\Theta\rangle}

\newcommand{\Sigmai}{\Sigma_{\rm i}}

\newcommand{\polya}{{\rm S}(\liecomplex{a})}

\newcommand{\ma}{m_\a}
\newcommand{\mduea}{m_{2\a}}

\newcommand{\pedtheta}[1]{{#1}_{\scriptscriptstyle{\Theta}}}

\newcommand{\Pavtheta}{\pedtheta{\rm P}^{\rm av}}

\newcommand{\ped}[2]{{#2}_{\scriptscriptstyle{#1}}}

\newcommand{\pedPi}[1]{{#1}_{\scriptscriptstyle{\Pi}}}

\newcommand{\cthetampl}{\pedtheta{c}^+(m;\l)}
\newcommand{\cthetamml}{\pedtheta{c}^-(m;\l)}
\newcommand{\cPi}{\pedPi{c}}
\newcommand{\phil}{\varphi_\l}

\newcommand{\varphiNCC}{{\varphi_{\scriptscriptstyle{\Pi_0}}(\lambda,a)}}
\newcommand{\cNCC}{{\ped{\Pi_0}{c}}}
\newcommand{\hyper}[4]{\ensuremath{\sideset{_{_2}}{_{_1}}{\mathop{F}}
\left(#1,#2;#3;#4\right)}}

\newcommand{\cCa}{C_c^\infty(C)}
\newcommand{\cCA}{C_c^\infty(C)}

\newcommand{\pwthetamC}{\pedtheta{\rm PW}(m;C)}

\newcommand{\ccitheta}{C_c^\infty(\pedtheta{A})^{\pedtheta{W}}}

\newcommand{\smC}{\scriptscriptstyle{\mathbb{C}}}

\begin{document}
\makeatletter
\title[Paley-Wiener theorems for the $\Theta$-spherical transform]
{Paley-Wiener theorems for the $\Theta$-spherical transform:\\ an overview}
\author{Gestur \'{O}lafsson}
\address{Department of Mathematics, Louisiana State University, Baton Rouge,
LA 70803, U.S.A.}
\email{olafsson@math.lsu.edu}
\thanks{The first author was partially supported by NSF grants DMS-0070607 and
DMS-0139783, and by the Sonderforschungsbereich Tr-12 
``Symmetrien und Universalit\"at in mesoskopischen Systemen''
of the Ruhr-Universit\"at Bochum.}
\thanks{Both authors were partially supported
by the Lorentz Center
at the Rijksuniversiteit Leiden.}
%
\author{Angela Pasquale}
\address{Fakult\"at f\"ur Mathematik, Ruhr-Universit\"at Bochum,
Universit\"atstrasse 150, Geb\"aude NA 4/65, 44780 Germany.}
\email{pasquale@cplx.ruhr-uni-bochum.de}
\date{}
\maketitle
\makeatother

\section{Introduction}

\noindent
Harmonic analysis has its origin in the work of Fourier on
the heat equation, which led him to consider the
expansion of an ``arbitrary'' 
$2\pi$-periodic function into superposition
of trigonometric functions:
$$f(x)\sim \sum_{n=-\infty}^\infty c_ne^{inx}$$
with
$$c_n=\frac{1}{2\pi }\int_{0}^{2\pi} f(t)e^{-int}\, dt\, .$$

One can interpret this expansion either as the spectral decomposition
of the differential operators with constant coefficients, or as decomposition of
$L^2([0,2\pi])$ into irreducible representations of the compact
Lie group $\mathbb{T}=\{z\in \mathbb{C} \mid  |z|=1\}$. There is
no ``definition'' of harmonic analysis that includes all
its different aspects, but the basic idea is \textit{to
study functions, or function spaces, in terms of
decomposition into simpler, basic functions}. This
includes spectral decomposition of differential
operators, theory of special functions,
integral transforms related to special functions,
atomic decomposition of function spaces, and
study of functions defined on a Lie group (or
a homogeneous space) by decomposing them into pieces associated with
the unitary irreducible representations of the group.
The meaning of ``simpler'' or ``basic'' functions depends then on the
context in which we are working.

The basic example of noncompact 
Lie group is the real line $\R$ considered as an additive group.
In this case all the unitary
irreducible representations are one dimensional and given
by the exponential functions $t \mapsto e^{i\l t}$ with $\l \in \R$.
To each (sufficiently regular) function $f:\R\to\C$
we associate its Fourier transform
\begin{equation*}
(\mathcal Ff)(\l)\sim \frac{1}{2\pi} \int_{\R} f(x) e^{-i\l x}\; dx.
\end{equation*}
Here the Fourier inversion
\begin{equation*}
f(x)\sim \int_{\R} (\mathcal Ff)(\l) e^{i\l x} d\l
\end{equation*}
provides the required decomposition with respect to the irreducible
representations or the spectral decomposition of the differential operators
on $\R$ with constant coefficients. 
Harmonic analysis also studies in which
sense this decomposition has to be considered. For instance,
the symbol ``$\sim$'' means convergence in $L^2$ norm if $f \in L^2(\R)$,
and even pointwise convergence if $f$ is smooth and compactly supported.

One of the fundamental questions
in harmonic analysis is to determine the
image of different function spaces under the Fourier transform.
For smooth, compactly supported functions, the
answer is given by the classical theorem of Paley
and Wiener, which characterizes this image in terms of holomorphic
extendibility and growth conditions.

The classical one-dimensional Fourier analysis on $\R$
has several generalizations.
The real line can be replaced by a higher-dimensional Euclidean space,
by a locally compact Hausdorff topological group, by a Lie group, or
by a homogeneous space. Among the homogeneous spaces, 
the symmetric spaces play an
important role for their numerous applications to other branches of
mathematics and to physics.

The aim of this paper is to give an overview of
several types of Paley-Wiener theorems, leading up
to the Paley-Wiener theorem for the $\Theta$-spherical Fourier transform.
The theory of $\Theta$-spherical functions is relatively new, and
originates from the interplay of the harmonic analysis on symmetric spaces
and the theory of special functions associated with root systems.
As of now, only some basic theorems, many of them only
for special cases, have been proven, leaving most
of the theory to be developed. The basic functions, the $\Theta$-spherical
functions, have singular behaviour as we approach the
boundary of the domain where they live.
Therefore the space of compactly supported functions
is easier to handle
than the $L^p$-spaces and many other natural function
spaces, and it is therefore
natural to look for a Paley-Wiener-type theorem.
Up to now such a theorem 
has been proven only in very special cases, and still,
its formulation and proof are very technical. In this paper we
shall not go into details of the proofs, but present an overview which
explains the different examples which have inspired and motivated
the theory of $\Theta$-spherical functions.
We refer to \cite{OP1}, \cite{OP2}, \cite{P1}, \cite{P2} and \cite{OP3} 
for details.

Starting from the Euclidean case,
we move to Harish-Chandra's theory of spherical functions on
Riemannian symmetric spaces of noncompact type.
The latter provided the geometric background to Heckman-Opdam's theory of
hypergeometric functions associated with root systems, which was developed
from the late 80ies. The theory of Heckman and Opdam not only yielded a
natural class of geometrically motivated special functions, but also
provided new methods and tools for understanding the geometric theory.

A generalization of Harish-Chandra's theory
to a different class of symmetric spaces found its roots in the
theoretical physics. In fact, the studies of the causal
structures of space-time underlined the physical relevance
of non-Riemannian symmetric spaces which are endowed with 
a $G$-invariant ordering. A special class of ordered symmetric spaces 
where a theory of spherical functions 
could be developed are the noncompactly causal (NCC) symmetric spaces.
This theory started in the early 90ies with the article \cite{FHO}.

It turned out that the Heckman-Opdam theory was the ideal tool to extend 
the results of \cite{FHO}.
At this point it seemed therefore natural to enclose all these different
but related theories in a single one, namely the theory of $\Theta$-spherical
functions.

The authors would like to thank G. van Dijk and V.F. Molchanov, 
the organizers
of the conference \textit{Representations of
Lie groups, harmonic analysis on homogeneous spaces
and quantization}, Leiden, December 9--13, 2002, and
the staff at the Lorentz Center, for hospitality during our
stay in Leiden.

\section{The Euclidean Paley-Wiener theorem}
\label{section:euclPW}

\noindent
Let $\frak a$ be an $r$-dimensional Euclidean real vector space, and
let $\liecomplex{a}$,
$\frak a^*$ and $\frakacs$ respectively denote the
complexification, the real dual and the
complex dual of the vector space $\frak a$.
The Euclidean Fourier transform of a sufficiently regular
function $f: \frak a \to \C$ is the function
$\mathcal Ff:\frakacs \to \C$ defined by
\begin{equation*}
 \mathcal Ff(\l):=\int_{\frak a} f(H)\, e^{\l(H)} \; dH,
\end{equation*}
where $dH$ denotes the Lebesgue measure on $\frak a$.
The integral converges for instance when $f$ is smooth and
compactly supported.

For a compact subset $E$ of $\frak a$ let $\conv(E)$ denote
the closed convex hull of $E$, i.e. the intersection of all closed half-spaces
in $\frak a$ containing $E$.
The support function of $E$ is the function $q_E:\frak a^* \to \R$
defined by
\begin{equation}
  \label{eq:support}
  q_E(\l):=\sup_{H \in E} \l(H)=\sup_{H \in \conv(E)} \l(H).
\end{equation}
Let $C$ be a compact convex subset of $\frak a$, and let $\cCa$ denote the
space of smooth functions on $\frak a$ with support contained in $C$.
The Paley-Wiener space ${\rm PW}(C)$ is the space of the entire functions
$g:\frakacs\to \C$
which are of exponential type $C$ and rapidly decreasing,
i.e., for every $N \in \N$
there is a constant $C_N\geq 0$ such that
$$
\abs{g(\l)} \leq C_N (1+\abs{\l})^{-N} e^{q_C(\Re\l)}
$$
for all $\l \in \frakacs$.

The classical theorem Paley and Wiener
characterizes ${\rm PW}(C)$ as
the image of $\cCa$ under the Euclidean Fourier transform
(see e.g. \cite{HoerLPDO}, Theorem 7.3.1, or \cite{JL},
Theorem 8.3 and Proposition 8.6).

\begin{Thm}[Paley-Wiener] \label{thm:classicalPW}
Let $C$ be a compact convex subset of $\frak a$.
Then the Euclidean Fourier transform
maps  $\cCa$ bijectively onto ${\rm PW}(C)$.
Moreover, if $C$ is stable under the action of a finite group $W$
of linear automorphisms of $\frak a$, then the Fourier transform
maps the subspace ${\cCa}^W$ of $W$-invariant elements in $\cCa$
onto the subspace ${\rm PW}(C)^W$ of $W$-invariant elements in ${\rm PW}(C)$.
\end{Thm}

Suppose that $A$ is a connected simply connected  abelian Lie  group with Lie algebra $\frak a$.
Then $\exp:\frak a\to A$ is a diffeomorphism. Denote by $\log$  the inverse of $\exp$.
The Euclidean Fourier transform of a sufficiently regular
functions $f: A \to \C$  is the function  $\mathcal F_Af: \frakacs \to \C$ defined by
\begin{equation} \label{eq:FA}
 \mathcal F_Af(\l):=\int_{A} f(a) e^{\l(\log a)} \; da,
\end{equation}
where the Haar measure $da$ on $A$ is the pullback under the exponential map of the
Haar measure $dH$ on $\frak a$.
Let $W$ be a finite group acting on $\frak a$ by linear automorphism.
Then we define an action by $W$ on $A$ by $w (\exp H)=\exp w(H)$.
Denote by $\cCA$  the space of smooth functions on $A$ with
support in the compact set $\exp C$, and let $\cCA^W$ be the subspace of $W$-invariant
elements. Composition with $\exp$ and Theorem \ref{thm:classicalPW} prove that
the Euclidean Fourier transform $\mathcal F_A$ is a bijection of
$\cCA$ onto ${\rm PW}(C)$ and restricts to a bijection of
$\cCa^W$ onto ${\rm PW}(C)^W$.

\section{Harmonic analysis on Riemannian symmetric spaces}
\label{section:Rssp}

\noindent
In this section we present the theory of spherical 
functions on Riemannian symmetric spaces as developed by Harish-Chandra,
Helgason, Gangolli and others, and state Helgason-Gangolli's 
Paley-Wiener theorem for the associated spherical Fourier transform.  
The main references are \cite{He2} and \cite{GV}.

Let $G$ be a connected noncompact semisimple Lie group with finite center,
and let $K$ be a maximal compact subgroup of $G$.
Then the quotient manifold $G/K$ can be endowed with the 
structure of Riemannian symmetric space of the noncompact type. 
Harish-Chandra's theory of spherical functions on
$G/K$  is the harmonic analysis of the $K$-invariant
functions on $G/K$.  
(For technical reasons, the theory was in fact developed in 
the wider setting of
$G$ of the Harish-Chandra class. See \cite{GV}, Chapter 2.) 

There are several equivalent ways
to define \textit{spherical functions} on $G/K$. See e.g. \cite{He2}, Ch.~IV, \S\S~2--3.
For this, let $dk$ be a normalized Haar measure on
$K$ and let $\mathbb{D}(G/K)$ denote the (commutative) algebra
of invariant differential operators on $G/K$. We
identify functions on $G/K$ with right $K$-invariant functions
on $G$.
\medskip

\noindent
\textit{The integral equations:\;} Let $\psi$ be a complex-valued continuous 
function on $G$, not identically zero. Then $\psi$ is spherical if
\begin{equation}\label{eq-int}
\int_K\psi (xky)\, dk=\psi(x)\psi(y)
\end{equation}
for all $x,y\in G$.

Notice that the integral equations (\ref{eq-int}) 
ensure that $\psi$ is $K$-biinvariant and that $\psi(e)=1$,
where $e$ denotes the unit element of $G$.
\medskip

\noindent
\textit{The differential equations:\;} A left $K$-invariant complex-valued smooth function
$\psi$ on $G/K$ is spherical if $\psi (eK)=1$ and if 
there exists a character $\chi :\mathbb{D}(G/K)\to
\mathbb{C}$ such that  for all $D\in \mathbb{D}(G/H)$ we have
\begin{equation} \label{eq-diffe}
D\psi =\chi (D)\psi\, .
\end{equation}

Thus the spherical
functions are the normalized joint eigenfunctions of $\mathbb{D}(G/K)$. 
Notice that, since
$\D(G/K)$ contains the Laplace operator, which is elliptic,
all spherical functions are indeed real analytic functions on $G$.
\medskip

\noindent
\textit{Character characterization:\;} A left $K$-biinvariant complex-valued function
$\psi$ on $G$ is spherical if the map
\begin{equation}\label{eq-spfct}
f\mapsto \hat{f}(\psi):=\int_{G/K}f(x)\psi (x)\, dx
\end{equation}
is a homomorphism of the (commutative) convolution algebra $C_c^\infty(G/K)^K$ 
of $K$-biinvariant complex-valued functions on $G$.
\medskip

In short, the spherical functions play for the harmonic analysis of
radial (that is, $K$-invariant) functions on $G/K$ the same role as 
the exponential functions for the harmonic analysis on the real line.

Let $\theta :G \to G$ be the Cartan involution on $G$ corresponding to
$K$, i.e., $K=\{k\in G \mid \theta (k)=k\}$. Denote by the same
letter the derived involution $\theta :\mathfrak{g}\to \mathfrak{g}$.
Then $\fg=\fk\oplus \fp$ where
$$\fk=\{X\in\fg\mid \theta(X)=X\}$$
is the Lie algebra of $K$ and
$$\fp=\{X\in \fg\mid \theta{X}=-X\}\, .$$
Let $\fa$ be a maximal abelian subspace of $\fp$. Let $\Sigma =\Sigma (\fg,\fa)$
be the set of roots of $\fa$ in $\fg$. For $\a\in \Sigma$ let
$\fg_\a :=\{X\in \fg\mid \text{$[H,X]=\a (H)X$ for all $H\in\fa$}\,\}$ be
the corresponding root space, and set
$\ma :=\dim\fg_\a$. Let
$\fm=\{X\in\fk\mid [\fa ,X]=\{0\}\}$ be the centralizer of
$\fa $ in $\fk$. Then
$$\fg =\fm \oplus \fa \oplus\bigoplus_{\a\in \Sigma}\fg_\a\, .$$

The set $\{H\in \fa \mid \text{$\a (H)=0$ for some $\a \in \Sigma$\,}\}$
is a finite union of hyperplanes. We can therefore choose an
$X\in \fa$ such that $\alpha (X)\not= 0$ for all $\a \in \Sigma$.
Let $\Sigma^+:=\{\a \in \Sigma \mid \a (X)>0\}$. Then $\Sigma^+$
is a \textit{positive system} of roots. Let $\Pi$ be
the corresponding fundamental system of \textit{simple roots}. Let
$$\fa^+:=\{H\in \fa\mid \text{$\a(H)>0$ for all $\a \in \Sigma^+$}\}$$
and
$$\fn:=\bigoplus_{\a \in \Sigma^+}\fg_\a\, .$$
Then $\fn$ is a nilpotent Lie algebra. We set
$A:=\exp \fa$, $A^+:=\exp \fa^+$ and $N:=\exp (\fn)$. The centralizer and the  
normalizer of $A$ in $K$ are respectively $M=Z_K(A)$ and $M^\prime =N_K(A)$.

\begin{Thm}\label{eq-Iwasawa}
The map
\begin{equation}
K\times A\times N\ni (k,a,n)\mapsto kan\in G
\end{equation}
is an analytic diffeomorphism.
\end{Thm}

The decomposition of Theorem 
 \ref{eq-Iwasawa} is known as \textit{Iwasawa decomposition} of $G$.
Define $a_K : G\to A$ by
\begin{equation}\label{eq-IwasawaA}
x\in Ka_K(x)N
\end{equation}

\begin{Thm}\label{th-Cartandec} We have $G=KAK$ and
the map
$K/M\times A^+ \ni (kM,a)\mapsto ka\in G/K$
is an analytic diffeomorphism onto an open
dense subset of $G/K$.
\end{Thm}

The Weyl group $W:=M^\prime /M$ is a finite reflection group acting on $\fa$
and  --~by duality~-- on $\fa^*$.
Fix a $W$ invariant inner product $\inner{\cdot}{\cdot}$ on $\fa$. For
$\a \in \Sigma$ define $H_\a \in [\fg_\alpha,\fg_{-\a}]$ by
$\a (H_\a)=2$. Let $r_\a :\fa^* \to \fa^*$ be the reflection
$r_\a (\lambda )=\lambda -\lambda (H_\a )\a$. Then $W$ is
generated by $\{r_\a \mid \a\in \Pi\}$.   For $a=\exp (H)\in A$ and
$\lambda \in \fa_\C^*$ let
$$a^\lambda :=e^{\l (H)}$$
and notice that the homomorphism $A\to \C^*$ are exactly the
maps $a\mapsto a^\l$. 

The complexification and complex dual of $\frak a$ are respectively denoted by 
$\liecomplex{a}$ and $\frakacs$. We extend the inner product of $\frak a$ to 
$\frak a^*$ by duality and to $\liecomplex{a}$ and $\frakacs$ by $\C$-bilinearity.
Finally, we set 
\begin{equation}
  \label{eq:rho}
  \rho=\frac{1}{2} \sum_{\a \in \Sigma^+} m_\a \a\,.
\end{equation}

We have now set up all the notations to describe Harish-Chandra's results
on spherical functions.

\begin{Thm}[Harish-Chandra]\label{thm:HCsph}
Let $\l \in   \fa_{\C}^*$. Then the function
$\varphi_\lambda :G\to \C$ given by
\begin{equation}\label{eq-intrepspf}
\varphi_\l (x):=\int_K a(xk)^{\lambda-\rho}\, dk
\end{equation}
is a spherical function. It is a real analytic function of $x \in G$
and a holomorphic function of $\l \in \frakacs$. Every spherical function
on $G$ is of the form $\varphi_\lambda$ for some
$\lambda\in \fa_\C^*$. Furthermore
$\varphi_\lambda=\varphi_\mu$ if and only if there
exists $w\in W$ such that $\lambda=w\mu$.
\end{Thm}

Notice that, since $\frak a^+$ is the interior of a fundamental domain of $W$,
Theorem \ref{th-Cartandec} implies 
that a continuous $K$-biinvariant functions $f$ on
$G$ is uniquely determined by its restriction $\mathrm{Res}_{A^+}(f)=f|_{A^+}$
to $A^+$. Moreover, 
we can normalize the invariant measures $dx$ on $G/K$ 
and $da$ on $A$ so, that for $f\in C_c(G/K)^K$ we have
$$\int_{G/K} f(x)\; dx=\int_{A^+} f(a)\delta (a)\; da
=\frac{1}{\abs{W}} \int_{A} f(a)\delta (a)\; da\,,$$
where $\abs{W}$ denotes the cardinality of the Weyl group and 
\begin{equation} \label{eq:delta}
\delta (a):=  \prod_{\alpha\in\Sigma^+}\abs{a^\a-a^{-a}}^{m_\a}.
\end{equation} 
Furthermore, a $K$-biinvariant function $\varphi$
is an eigenfunction for $\mathbb{D}(G/K)$ if and only if
$\mathrm{Res}_{A^+}(f)$ is an eigenfunction for
the system of equation on $A^+$ given by the radial components
of operators from $\mathbb{D}(G/K)$.
Thus harmonic analysis of $K$-invariant functions on
$G/K$ becomes the study of $W$-invariant functions on $A$ or even of 
functions on the set $A^+$. Let $\D_{G/K}(A)$ denote the commutative 
algebra of radial components along $A^+$ of the differential operators in
$\D(G/K)$. With each $\l \in \frakacs$ is associated a
character $\chi_\l$ of $\D_{G/K}(A)$, and the spherical function
$\varphi_\l$ is determined on $A^+$ by the system of partial differential equations
\begin{equation}
  \label{eq:diffsphHC}
  D\varphi=\chi_\l(D)\; \varphi, \qquad D \in \D_{G/K}(A).
\end{equation}
 Observe that $\chi_\l=\chi_{w\l}$ for all $w\in W$.

Let $\{H_i\}_{i=1}^r$ be a fixed orthonormal basis of $\frak a$. For each 
$H \in \frak a$ let $\partial(H)$ denote the corresponding 
directional derivative in $\frak a$.  Then the
partial differential equation (\ref{eq:diffsphHC}) on $A^+$ corresponding to 
the Laplace operator can be written explicitly as
\begin{equation}\label{eq:diffphil}
L\varphi=(\inner{\l}{\l}-\inner{\rho}{\rho}) \varphi
\end{equation}
where
\begin{equation}
  \label{eq:radiallaplace}
  L:=\sum_{i=1}^r \partial(H_i)^2+ \sum_{\a \in \Sigma^+} m_\a 
\frac{1+e^{-2\a}}{1-e^{-2\a}} \; \partial(A_\a)
\end{equation}
denotes the radial component on $A^+$ of the Laplace operator.
A multidimensional variant of the classical method of Frobenius for 
determining local solutions of differential equations with regular 
singularities led Harish-Chandra to look for solutions of (\ref{eq:diffphil})
of the form 
\begin{equation}\label{eq:HCseries}
\Phi_\l(m;a)=a^{\l-\rho} 
\sum_{\mu \in \Lambda} \Gamma_\mu(m;\l) a^{-\mu}, 
\qquad a \in A^+.
\end{equation}
Here $\Lambda:=\left\{\sum_{j=1}^r n_j \a_j\mid n_j \in \N_0 \right\}$, where
 $\a_1,\dots, \a_r$ is for some
enumeration of $\Pi$.
Moreover $m:=\{m_\a\mid\a \in \Sigma\}$ denotes the set of multiplicities.  
With the initial condition $\Gamma_0(m;\l)=1$, 
the coefficients $\Gamma_\mu(m;\l)$ with $\mu \in \Lambda \setminus \{0\}$ 
are uniquely determined by means of the recurrence relations
\begin{equation} \label{eq:recursion}
\inner{\mu}{\mu-2\l} \Gamma_\mu(m;\l)= 2 \sum_{\a\in \Sigma^+} m_\a 
 \sum_{\substack{k\in \N\\\mu -2k\a \in \Lambda}}
 \Gamma_{\mu-2k\a}(m;\l)  \inner{\mu+\rho-2k\a-\l}{\a}
\end{equation}
provided $\l \in \frakacs$ satisfies $\inner{\mu}{\mu-2\l} \neq 0$ for all 
$\mu \in  \Lambda \setminus \{0\}$. On this subset of $\frakacs$ the series
on the right-hand side of (\ref{eq:HCseries}) converges 
to a meromorphic function of $\l$ and 
real analytic function of $a \in A^+$ -- in fact, there is tubular neighborhood $U^+$ 
of $A^+$ in the comlpexification of $A$ on which this series converges to a holomorphic 
function --. \label{Phiholo}
It is remarkable that the function $\Phi_\l$
turns out to solve the entire system (\ref{eq:diffsphHC}) and allows to construct a 
basis for the smooth solutions of (\ref{eq:diffsphHC}) on $A^+$.

\begin{Thm}[Harish-Chandra]\label{th34}
Let the notation be as above. Then the following properties hold
for generic spectral parameter $\lambda\in \fa_\C^*$.
\begin{enumerate}
\thmlist
\item Let $D\in \D_{G/K}(A)$. Then $D\Phi_{\lambda}=\chi_\lambda (D)\Phi_\l $.
\item  The functions $\{\Phi_{w\lambda}\mid w\in W\}$
form a basis for the space of smooth solutions to the
system (\ref{eq:diffsphHC}) of differential equations on $A^+$.
\item 
There is a meromorphic function $c$ on $\frakacs$ (depending only on $\l$ and 
on the root structure) 
so that the 
spherical function   
$\phil$ admits on $A^+$ the expansion
\begin{equation}
\label{eq:orHCexpansion}
 \phil(a)=\sum_{w \in W} c(w\l) \Phi_{w\l}(m;a)\,.
\end{equation}
\end{enumerate}
\end{Thm}

The function $c$ occurring in (\ref{eq:orHCexpansion}) is the so-called Harish-Chandra's 
$c$-function. It governs the asymptotic behavior of $\phil$ on $A^+$. More precisely,
let us write $A^+\ni a \to \infty$ to indicate that $a\in A^+$ and
$\lim a^{-\a }=0$ for all $\a\in\Sigma^+$. If $\Re\inner{\lambda +\rho}{\alpha} <0$ for
all $\alpha \in \Sigma^+$, then
$$
\lim_{A^+\ni a  \to \infty}  a^{\rho -\lambda}\varphi_\lambda (a)=c(\lambda)\, .$$
The $c$-function is given by 
\begin{equation*}
  c(\l)=\int_{\bar{N}} e^{-(\l+\rho)(H(\bar{n}))}\; d\bar{n},
\end{equation*}
where $\bar{N}:=\theta(N)$ and the Haar measure $d\bar{n}$ on $\bar{N}$
is normalized by the condition 
\begin{equation} \label{eq:crho}
  c(\rho):=\int_{\bar{N}} e^{-\rho(H(\bar{n}))}\; d\bar{n}=1\,.
\end{equation}
Let $\Sigmai^+$ denote the set of indivisible positive roots. 
Hence $\a \in \Sigmai^+$
provided $\a \in \Sigma^+$ but $\a/2\notin \Sigma^+$.
For $\l \in \frakacs$ set 
\begin{equation}
  \label{eq:la}
  \l_\a:=\frac{\inner{\l}{\a}}{\inner{\a}{\a}}.
\end{equation}
Then Gindikin and Karpelevic proved the following explicit product 
formula for Harish-Chandra's 
$c$-function:
\begin{equation}
  \label{eq:HCc}
  c(\l)= \kappa_0 \prod_{\a \in \Sigmai^+}
\frac{2^{-\l_\a} \; \Gamma\left(\l_\a \right)}
{\Gamma\left( \frac{1}{2}\big(\l_\a+m_\a/2+1\big)\right)
\Gamma\left(\frac{1}{2}\big(
\l_\a +m_\a/2+m_{2\a}\big)\right)},
\end{equation}
where the constant $\kappa_0$ is chosen so that (\ref{eq:crho}) holds.
See e.g. \cite{GV}, Theorem 4.7.5,
or \cite{He2}, Ch. IV, Theorem 6.14. 

\begin{Ex}[The real rank-one case]
\label{ex:rankone}
The real rank-one case corresponds to 
Riemannian symmetric spaces of the noncompact type 
for which $\frak a$ is one dimensional. The set $\Sigma^+$ consists
at most of two elements: $\a$ and, possibly, $2\a$. 
By setting $H_\a/2\equiv 1$ and $\a\equiv 1$,
we identify $\frak a$ and $\frak a^*$ with $\R$, and  their complexifications
$\liecomplex{a}$ and $\frakacs$ with $\C$. 
The Weyl chamber $\frak a^+$ coincides with the half-line $(0,+\infty)$. 
The Weyl group $W$ reduces to $\{-1,1\}$ acting on $\R$ and $\C$ by multiplication. 
Moreover  $\rho(m)\equiv \ma/2+\mduea$. We normalize the inner product so that 
$\inner{\a}{\a}=1$. The algebra $\D(G/K)$ is generated by the Laplace operator, and
the system of differential equation (\ref{eq:diffphil}) with spectral parameter $\l \in \C$ 
is equivalent to the single Jacobi differential equation
\begin{equation} \label{eq:hyperalge}
\frac{d^2 \varphi}{dt^2}+\big(\ma\coth t+\mduea\coth(2t)\big)\; \frac{d\varphi}{dt} 
=(\l^2-\rho^2)\varphi, \qquad t \in (0,+\infty).
\end{equation}
The solution to (\ref{eq:hyperalge}) that behaves asymptotically as $e^{(\l-\rho)t}$ 
for $t\rightarrow +\infty$ is 
\begin{equation*}
  \label{eq:Phillrankone}
\Phi_\l(m;t)
=(2\sinh t)^{\l-\rho} \hyper{\frac{\rho-\l}{2}}{\frac{-\ma/2+1-\l}{2}}{1-\l}{-\sinh^{-2}t}, 
\end{equation*}
where $\sideset{_{_2}}{_{_1}}{\mathop{F}}$ denotes the Gaussian hypergeometric function.
The function $\Phi_\l(m;t)$ coincides with the Jacobi function of second kind
$\Phi^{(a,b)}_\nu(t)$ with parameters 
$a=(\ma+\mduea-1)/2$, $b=(\mduea-1)/2$ and $\nu=-i\l$. See e.g. \cite{Koo84}, 
Section 2.
Harish-Chandra's spherical functions are Jacobi functions of first kind
\begin{equation*}
  \varphi_\l(t)=
\hyper{{\frac{\ma+2\mduea+2\l}{4}}}{{\frac{\ma+2\mduea-2\l}{4}}}{{\frac{\ma+2\mduea+1}{2}}}
{-\sinh^2 t},
\end{equation*}
and the expansion (\ref{eq:orHCexpansion}) reduces to 
one of Kummer's relations between  
solutions of the hypergeometric equations. See e.g. \cite{Er}, 2.9 (34).
\end{Ex}

\begin{Ex}[The complex case] \label{ex:complex}
The complex case corresponds to Riemannian symmetric spaces $G/K$ 
for which the Lie algebra $\frak g$ of $G$ 
is endowed with a complex structure. It is 
characterized by the fact that the root system $\Sigma$ is reduced
with all multiplicities $m_\a=2$. It follows, in particular, that 
$\rho=\sum_{\a\in\Sigma^+} \a$.  
The Harish-Chandra series are given in $A^+$ by 
$
  \Phi_\l(2;a):=\Delta(a)^{-1}a^{\l}
$
where 
\begin{equation}\label{eq:Weylden}
\Delta(a):=\prod_{\a\in \Sigma^+} (a^\a-a^{-\a})
\end{equation}
is Weyl's denominator. Set
\begin{equation}
\label{eq:pi}
\pi(\l):=\prod_{\a\in\Sigma^+} \l_\a.
\end{equation} 
An explicit formula for the spherical functions on $G/K$ was determined by 
Harish-Chandra, namely
\begin{equation}
  \label{eq:varphicompl}
  \varphi_\l(a)=\frac{\pi(\rho)}{\pi(\l)}\; 
\frac{\sum_{w\in W} (\det w)
 a^{w\l}}{\Delta(a)}
\end{equation}
(see e.g. \cite{He1}, Ch. IV, Theorem 5.7).
Notice that $c(\l)=\pi(\rho)/\pi(\l)$.
\end{Ex}

The \textit{spherical Fourier transform} $\mathcal Ff$ of 
$f \in C_c(G/K)^K$ is  
the $W$-invariant function on $\frakacs$ defined by
\begin{equation}
  \label{eq:sphtr}
  (\mathcal F f)(\l):=\int_G  f(x)\phil(x) \; dx =
   \int_{A^+} f(a)\varphi_\l(a)\delta(a) \;da=\frac{1}{\abs{W}} \; \int_{A} 
f(a)\varphi_\l(a)\delta(a) \;da\,, 
\end{equation}
where $dx$ and $da$ are respectively 
the fixed normalizations of the Haar measures on $G$ and $A$, and $\delta$ is given by 
(\ref{eq:delta}).

\begin{Thm}[Plancherel Theorem] \label{th-Plan2}
The spherical Fourier transform extends to an isometric isomorphism
$$\mathcal{F}:L^2(G/K)^K\simeq L^2(A,\delta(a)da)^W 
\longrightarrow  L^2(i\fa,|c(\lambda)|^{-2}d\lambda)^W$$
with pointwise inversion on $C_c^\infty(G/K)^K$ given by
\begin{equation} \label{eq:invsphHC}
f(a)=\int_{i\fa^*}\mathcal{F}(f)(\lambda )\varphi_{-\lambda}(a)|c(\lambda)|^{-2}\, 
d\l\, ,\quad
a\in A\,.
\end{equation}
\end{Thm}

In terms of representation theory, this corresponds to
the decomposition of $L^2(G/K)$ into
direct integral of spherical principal series representations.

Recall from Section \ref{section:euclPW} the definition of Paley-Wiener space 
${\PW(C)}^W$ associated with a compact convex $W$-invariant subset $C$ of a 
Euclidean space. Here and in the following sections, 
$W$ always denotes the Weyl group. 

The characterization of the image of $C_c^\infty(G/K)^K$ 
is the content of the Paley-Wiener theorem for the spherical Fourier transform. 
As in the Euclidean case, this theorem gives a finer characterization 
by describing, for every compact convex $W$-invariant subset 
$C$ of $\frak a$, the image of the space $C_c^\infty(C)^W$ of 
elements $f \in C_c^\infty(G/K)^K$ for which the support 
$\supp f|_A$ of their $W$-invariant restriction to $A$ is contained in $\exp(C)$. 
The Paley-Wiener theorem was stated by Helgason in 1966 in \cite{HePW}, 
where it was proven for the real rank-one and complex cases. 
The extension to arbitrary Riemannian symmetric spaces of the noncompact type was 
completed by Gangolli \cite{G1} in 1971. The proof was later simplified by 
Rosenberg \cite{Ros}. 

\begin{Thm}[Helgason-Gangolli-Rosenberg]
\label{thm:RiemPW}
Let $C$ be a $W$-invariant compact convex subset of $\frak a$. 
Then the Euclidean Fourier transform maps $C_c^\infty(C)^W$ 
bijectively onto ${\rm PW}(C)^W$.
\end{Thm}

The original formulation of Theorem \ref{thm:RiemPW} 
considers only the case in which $C$ is some
Euclidean ball $B_R:=\{H \in \frak a
\mid \abs{H}:=\inner{H}{H}^{1/2} \leq R\}$ with $R>0$. 
Its extension to
arbitrary $W$-invariant convex compact subsets of $\frak a$ is elementary.
See e.g. \cite{HS} or \cite{OP3}. 
Nevertheless, to simplify our exposition, here we only outline the proving
method for the case of $C=B_R$. 

The fact that the spherical Fourier transform maps
$C_c^\infty(B_R)^W$ into ${\rm PW}(B_R)^W$ 
is obtained by writing this transform as 
composition $\mathcal F=\mathcal F_A \circ \mathcal A$ of 
the so-called Abel transform 
and of the Euclidean Fourier transform (\ref{eq:FA}). 
Since the Abel transform can be easily seen to 
map $C_c^\infty(B_R)^W$ into itself, 
the required inclusion follows then from the Euclidean 
Paley-Wiener theorem (Theorem \ref{thm:classicalPW}). 
The thrust of the Paley-Wiener theorem for $\mathcal F$ 
is to get the surjectivity, namely,
that given a holomorphic function $h \in {\rm PW}(B_R)^W$, there is an 
$f \in C_c^\infty(B_R)^W$ so that $\mathcal Ff=h$. 
Motivated by the inversion formula (\ref{eq:invsphHC}), 
one defines $f$ to be a suitable constant multiple of the wave packet 
\begin{equation}
  \label{eq:IHC}
  \mathcal I h(x):=\int_{i\frak a^*} h(\l) \varphi_{-\l}(x) 
\abs{c(\l)}^{-2}\; d\l\,, \qquad x \in G.
\end{equation}
Then $\mathcal I h$ is well-defined, smooth and $K$-biinvariant.
The crucial step is to show that $\mathcal I h(\exp H)=0$ for $\abs{H}>R$. 
For this one
uses the expansion (\ref{eq:orHCexpansion}) 
and the $W$-invariance of $h$ to write
\begin{equation*}
 \mathcal I h(\exp H)=
\int_{i\frak a^*} h(-\l) c(-\l)^{-1} \Phi_\l(m;a) \; d\l\,.
\end{equation*}
At this step one would like to replace $\Phi_\l(m;a)$ by its definition (\ref{eq:HCseries}) 
and interchange summation and integration to get
\begin{equation} \label{eq:interchange}
  \mathcal I h(\exp H)=\sum_{\mu \in \Lambda} e^{-(\mu+\rho)(H)} B_\mu(m;H)
\end{equation}
with
\begin{equation}  \label{eq:tobeshifted}
  B_\mu(m;H):=\int_{i\frak a^*} h(-\l) 
\Gamma_\mu(m;\l) c(-\l)^{-1} e^{\l(H)} \; d\l\,.
\end{equation}
Now the recursion relations (\ref{eq:recursion}) and the product 
formula (\ref{eq:HCc}) for $c(\l)$ show that both 
$\Gamma_\mu(m;\l)$ and $c(-\l)^{-1}$ are holomorphic in $\l$ provided 
$\Re\inner{\l}{\a}<0$ for all $\a \in\Sigma^+$. This, together with suitable estimates
for $\Gamma_\mu(m;\l)$, allows one to shift the contour of integration in (\ref{eq:tobeshifted})
and to prove, as in the Euclidean Paley-Wiener theorem, that $B_\mu(m;H)=0$ 
for $\abs{H}>R$. 
The estimates that have been essential to justify (\ref{eq:interchange}) were established in
\cite{G1}.

\section{The Heckman-Opdam Theory}  \label{section:HO}

\noindent
As noticed in Example \ref{ex:rankone}, 
the spherical functions on rank-one 
Riemannian symmetric spaces are special instances of 
Jacobi functions of first kind, 
hence of hypergeometric functions. 
The specialization occurs with the choice 
of the multiplicities $m_\a$ and $m_{2\a}$ as certain 
nonnegative integers 
fixed by the geometry. 
Also in the higher-rank case, the geometry constraints 
the root multiplicities $m_\a$ to assume certain specific values.
The spherical functions as well are determined by the geometry,
since the system of differential equations (\ref{eq:diffsphHC})
originates from the algebra of $G$-invariant
differential operators on $G/K$.
Nevertheless, the differential equation (\ref{eq:diffphil}) 
makes perfectly sense without the
geometrical restrictions on $m$. 
This observation was the starting point of the
theory of hypergeometric functions associated with root systems, 
which has been developed
by Heckman and Opdam.
Their goal was to reconstruct, for arbitrary 
complex values of multiplicities, the systems of differential equations 
(\ref{eq:diffsphHC}). 
As analytic continuations (in the multiplicity parameters) 
of Harish-Chandra's spherical 
functions, the common eigenfunctions of these new systems would have 
provided a class of 
geometrically motivated multivariable hypergeometric functions.
Heckman and Opdam could realize their program in a series of papers
from 1988 to 1995 (\cite{HOpd1},\cite{Heck1}, \cite{Opd1},\cite{Opd2}, 
\cite{Heck2}, \cite{OpdGauss}, \cite{OpdActa}). 
The generalized spherical functions they constructed 
are nowadays known as \textit{hypergeometric functions associated with 
root systems}.    
Our overview below is based mainly on \cite{HS}, 
\cite{HeckBou} and \cite{OpdActa}, to 
which we refer for details. 

The Riemannian symmetric spaces of Harish-Chandra's theory are replaced in the 
theory of Heckman and Opdam by triples $(\frak a,\Sigma,m)$, where 
$\frak a$ is an $r$-dimensional real Euclidean vector space, 
$\Sigma$ is a root system in 
$\frak a^*$, and $m$ is a multiplicity function on $\Sigma$,
 that is a function $m:\Sigma\rightarrow \C$ 
which is invariant under the Weyl group $W$ 
of $\Sigma$. Setting $m_\a:=m(\a)$ for $\a \in \Sigma$, we therefore have 
$m_{w\a}=m_\a$ for all $w \in W$.  Because of our interest 
in Paley-Wiener theorems, we
shall restrict ourselves here to the case in which all 
$m_\a$ are nonnegative reals.
Harmonic analysis results for some negative values of 
$m_\a$ can be found in \cite{Opd00}.

We say that the triple is \textit{geometric} if there is a 
Riemannian symmetric space of noncompact type $G/K$ 
with restricted root system $\Sigma$ for the 
corresponding pair $(\frak g,\frak a)$ 
so that $m_\a$ is the multiplicity of the root $\a$ for all $\a\in \Sigma$. 
Notice that we adopt the multiplicity notation 
commonly used in the theory of symmetric spaces. 
It differs from the notation employed by Heckman and Opdam in the 
following ways.
The root system $R$ used by Heckman and Opdam
is related to our root system $\Sigma$ by the relation 
$R=\{2\a\mid \a \in \Sigma\}$; 
the multiplicity function $k$ in Heckman-Opdam's work is related to 
our $m$ by $k_{2\a}=\ma/2$.

In the following our fixed inner product in $\frak a$ is denoted by 
$\inner{\cdot}{\cdot}$.  The dimension $r$ of $\frak a$ is called the 
(real) rank of the triple $(\frak a, \Sigma,m)$.
The symbols $\liecomplex{a}$, $\frakacs$, 
$H_\l$, $\frak a^+$, $\Sigma^+$, $\Pi$, $r_\a$, $\l_\a$ shall 
have the same meaning as in Section \ref{section:Rssp}.  
To construct an analog of the Cartan subgroup, 
we first consider the complex torus 
$\complex{A}:=\liecomplex{a}/\Z\{i\pi H_\a: \a \in \Sigma\}$ with Lie
algebra $\liecomplex{a}$. 
Let $\exp:\liecomplex{a} \to \complex{A}$ denote the canonical
projection. 
Then $A:=\exp \frak a$ is an abelian group so that $\exp: \frak a \to A$ is a 
diffeomorphism. 
We write $\log$ for the inverse of $\exp$, and set $A^+:=\exp \frak a^+$.

The restricted weight lattice of $\Sigma$ is the set $P$ of 
all $\l\in \frak a^*$ for which $\l_\a \in \Z$ for all $\a\in \Sigma$. 
Notice that $2\a \in P$ for all $\a \in \Sigma$ and that 
$P$ consists exactly of the
elements $\l \in \frakacs$ for which the exponential 
$e^\l(a):=e^{\l(\log a)}=a^\l$ is single valued 
on $\complex{A}$.  We denote by $\C[\complex{A}]$ the
$\C$-linear span of the $e^\l$ with $\l\in P$.

Let $\polya$ denote the symmetric algebra over $\liecomplex{a}$ 
considered as the space of polynomial functions on $\frakacs$, 
and let $\polya^W$ be the subalgebra of $W$-invariant elements.
For $p\in \polya$ write $\partial (p)$ for the corresponding
 constant-coefficient differential operator
on $\complex{A}$ (or on $\liecomplex{a}$). 

An algebra of differential operators playing the role of 
the algebra $\D_{G/K}(A)$ was obtained by Heckman and Opdam as follows. 
Let $\mathcal R$ 
be the subalgebra of the quotient field of $\C[\complex{A}]$ 
generated the constant function 
$1$ and by $1/(1-e^{-2\a})$ with $\a \in \Sigma^+$. Let 
$\D\mathcal R^W$ denote the algebra of 
differential operators on $\complex{A}$ with coefficients 
in $\mathcal R$ which are invariant under the Weyl group $W$.
If $L(m)$ denotes the differential operator defined 
(for arbitrary multiplicity functions
$m$) by the right-hand side of (\ref{eq:radiallaplace}),  
then $L(m) \in \D\mathcal R^W$.
The elements in $\D\mathcal R^W$ which commute with $L(m)$ form therefore 
an algebra $\D(\frak a, \Sigma,m)$ which agrees with  
$\D_{G/K}(A)$ when $(\frak a, \Sigma,m)$ is the geometric triple which
corresponds to the Riemannian symmetric space 
of noncompact type $G/K$. 
Heckman and Opdam proved that $\D(\frak a, \Sigma,m)$ 
is indeed commutative and parameterized by the elements of $\polya^W$. 
It follows, in particular, that is 
$\D(\frak a, \Sigma,m)$ generated by $r(=\dim \frak a)$ 
elements. 
It is important to point out that the elements 
of $\D(\frak a, \Sigma,m)$ can
be constructed algebraically. 
Indeed Cherednik determined an algebraic algorithm for 
constructing the operator $D(m;p) \in \D(\frak a, \Sigma,m)$ 
corresponding to the 
element $p \in \polya^W$ directly from $p$ and the data $(\frak a, \Sigma,m)$. His main tool are the so-called Dunkl-Cherednik 
differential-reflection operators.
The operators $D(m;p)$ are known as the 
\textit{hypergeometric differential operators}.   

Let  $\l \in \frakacs$ be arbitrarily fixed. 
Then the system of differential equations
\begin{equation}
  \label{eq:hypereq}
  D(m;p) \varphi=p(\l)\varphi, \qquad p \in \polya^W
\end{equation}
is called the \textit{hypergeometric system of differential equations 
with spectral parameter $\l$} associated with the data $(\frak a,\Sigma,m)$.
For geometric triples it agrees
with the system of partial differential equations (\ref{eq:diffsphHC}) defining 
Harish-Chandra's spherical function
of spectral parameter $\l$.

Formula (\ref{eq:HCseries}) can still be employed to define 
the Harish-Chandra series. 
We extend the definition of Harish-Chandra's $c$-function to arbitrary
multiplicity functions $m$ by means of the Gindikin-Karpelevic 
formula (\ref {eq:HCc}). 
To underline its dependence on $m$, we replace the notation 
$c(\l)$ by $c(m;\l)$. Similarly, we write $\delta(m;a)$ instead of 
$\delta(a)$ in (\ref{eq:delta}), and $\rho(m)$ instead of $\rho$ in
(\ref{eq:rho}).
   
Notice that $c(m;\l)$ is a meromorphic function of $\l \in \frakacs$. 
Also the Harish-Chandra series $\Phi_\l(m;a)$
are meromorphic in $\l \in \frakacs$ and singular on the walls of $A^+$.
So a priori the sum (\ref{eq:orHCexpansion}) is only well-defined 
on $A^+$ and meromorphic 
in $\l$. This singular behavior appears of course also in the 
Riemannian case, but, there,
the right-hand side of  (\ref{eq:orHCexpansion}) 
is known to be regular by Theorem 
\ref{thm:HCsph}.
Heckman and Opdam, who used the right-hand side of 
 (\ref{eq:orHCexpansion}) to {\em define} their
spherical functions, needed to develop tools to understand how
the cancellation of singularities takes place.  They proved the following 
fundamental result (cf. e.g. Theorem 4.4.2 in \cite{HS}).

\begin{Thm}[Heckman and Opdam] \label{thm:HOHCexp} 
For a fixed nonnegative multiplicity function $m$ there is a 
$W$-invariant tubular neighborhood $U$ of $A$ in $\complex{A}$ so that the function 
\begin{equation*}
 \varphi_\l(m;a):=\sum_{w \in W} c(m;w\l) \Phi_{w\l}(m;a), \qquad a\in A^+,
\end{equation*}
extends to a $W$-invariant holomorphic function of 
$(\l,a) \in \frakacs \times U$.
\end{Thm}

The functions $\phil(m;a)$ are the so-called
\textit{hypergeometric functions associated with the triple $(\frak a,\Sigma,m)$} 
(or, simply, \textit{with root system $\Sigma$}). By construction, they reduce to
Harish-Chandra's spherical functions in the geometric case.

Let $C^\infty_c(A)^W$ denote the space of $W$-invariant $C^\infty$ functions
on $A$ with compact support. The \textit{Opdam transform} of $f \in C^\infty_c(A)^W$
associated with the triple $(\frak a,\Sigma,m)$ is the $W$-invariant function 
$\mathcal Ff(m)$ on $\frakacs$ defined by 
\begin{align*}
{\mathcal F}f(m;\l)&:=
\frac{1}{\abs{W}} \, \int_{A} 
f(a) \, \varphi_\l(m;a) \delta(m;a)\; da 
=\int_{A^+} f(a) \, \varphi_\l(m;a) \delta(m;a)\; da \,.
\end{align*}
Notice that the assumption $m_\a\geq 0$ for all $\a \in \Sigma$ 
ensures that $\delta(m;a)$ is a continuous $W$-invariant function of $a \in A$. 
The inversion formula for the Opdam transform is given by 
\begin{equation} \label{eq:inversionOpdam}
f(a)= \kappa 
\int_{i\frak a^*} {\mathcal F}f (m;\l) \, \varphi_{-\l}(m;a)  \abs{c(m;\l)}^{-2}\; d\l
\end{equation}
where $d\l$ is a suitable normalization of the Lebesgue measure on $i\frak a^*$ and
$\kappa >0$ is a constant depending only on  the normalization of the measures.

Recall from Section \ref{section:euclPW} the notation  
$C^\infty_c (C)^W$ for the
$W$-invariant $C^\infty$ functions 
$f:A\to\C$ with support contained in $\exp C$, where
$C$ is a compact, convex and $W$-invariant 
subset of $\frak a$. 
The Paley-Wiener theorem for the Opdam transform have been 
proven in \cite{OpdActa},
Theorems 8.6 and 9.13(4); see also \cite{Opd00}, p. 49.

\begin{Thm}[Opdam]
\label{thm:pwopdam}
Let $m$ be a fixed multiplicity function 
with $m_\a \geq 0$ for all $\a \in \Sigma$.
Suppose $C$ is a compact, convex and $W$-invariant 
subset of $\frak a$. Then ${\mathcal F(m)}$ maps 
$C^\infty_c (C)^W$ bijectively onto $\PW(C)^W$. 
\end{Thm}

In its original form, Theorem \ref{thm:pwopdam} 
was stated for $C=\conv W(H)$, where $W(H)$ is
the Weyl group orbit of an element $H \in \frak a$ and 
$\conv$ denotes the convex hull. The above formulation can be easily 
deduced from it.

It is important to observe that the proof of the surjectivity of 
Opdam transform follows Helgason's method. But the fact that 
${\mathcal F(m)}$ maps 
$C^\infty_c (C)^W$ into $\PW(C)^W$ presents new difficulties. The crucial tool
used in the Riemannian case was the decomposition 
$\mathcal F=\mathcal F_A\circ \mathcal A$. In the context of 
hypergeometric functions associated with root systems,
the Abel transform $\mathcal A$ does not have a direct extension
--~one could of course generalize the Abel transform as 
$\mathcal F^{-1}_A \circ \mathcal F(m)$, but this definition does not
help to determine the support of $\mathcal F(m)f$~--.
Opdam's method is based on estimates for the hypergeometric 
functions $\varphi_\l$. 
Compared to the Riemannian case,  the difficulty for estimating 
the hypergeometric functions relies essentially in the fact that the 
hypergeometric functions associated with root systems 
are generally not representable by integral formulas.
Opdam's estimates have been obtained by introducing new regular nonsymmetric 
eigenfunctions of the Dunkl-Cherednik operators with Weyl groups 
symmetrization equal to the $\varphi_\l$. 
The important feature is that the Dunkl-Cherednik operators 
associated with elements in 
$\liecomplex{a}$ are first order differential-reflection operators. 
Indeed Opdam was able to extend to this setting a clever method, due 
to de Jeu \cite{deJeu}, 
which uses the first order directional 
derivatives occurring in the eigenfunction problems to evaluate the 
growth of the eigenfunctions themselves.

\section{Spherical functions on noncompactly causal symmetric spaces}
\label{se-nonRsp}

\noindent
The theory of spherical functions on Riemannian symmetric spaces 
depends on at least three closely related facts. 
First, the subgroup $K$ is a maximal compact subgroup, i.e., $\theta$ is
a Cartan involution. Secondly, the algebra $\mathbb{D}(G/K)$
contains an elliptic differential operator, and hence all
the joint eigenfunctions -- or distributions -- are
real analyitc. Finally we have the decompositions
$G=KAN=KAK$. If we replace $\theta$ by an arbitrary
involution $\tau :G\to G$ and set $H=\{h\in G\mid \tau (h)=h\}$, then
$H$ is no longer compact, on $G/H$ there are in general no elliptic
invariant differential operators, $G\not= HAN$ (i.e., a corresponding
generalized Iwasawa decomposition does not hold in general),
and $G\not= HAH$. But it was shown in \cite{FHO} that there
is a natural class of symmetric spaces, the 
\textit{noncompactly causal symmetric
spaces}, where an analogue theory of spherical functions defined on
open $H$-invariant conal subset of $G/H$ can be developed.      
For simplicity we shall assume in the following that $G$ 
is contained in the simply connected complex Lie group $G_\C$ with
Lie algebra $\fg_\C=\fg\otimes_\R \C$.

Let $\tau: G\to G$ be a nontrivial involution commuting with the
fixed Cartan involution $\theta$. We assume that $\tau $ is not a
Cartan involution, so $H$ is not compact. As usual we denote the
differential of $\tau$ by the same letter. Let
$\fh :=\{X\in\fg\mid \tau (X)=X\}$ and 
$\fq:=\{X\in \fg\mid \tau (X)=-X\}$.
Then $\fh$ is the Lie algebra of $H$. We
have
$$\fg=\fh\oplus\fq=\fh_k\oplus\fh_p\oplus\fq_k\oplus \fq_p$$
where the subscript ${}_k$, respectively ${}_p$, denotes 
intersection with $\fk$, respectively $\fp$.
An element $X\in \fg$ is called \textit{hyperbolic}
if $\ad (X)$ is semisimple with real eigenvalues.

\begin{Def}\label{def-ncc} 
Assume that $G/H$ is simple. Then $G/H$ is
called \textit{noncompactly causal}, in short NCC, if there exists
an open convex $H$-invariant cone $\Omega\neq \emptyset$ 
in $\fq$, containing
no affine lines, such
that each $X\in \Omega$ is hyperbolic.
\end{Def}

We assume from now on that $G/H$ is NCC and that $\Omega$ is an
$H$-invariant cone as in Definition \ref{def-ncc}.
We shall assume that $\Omega$ is maximal.
We refer to \cite{HObook} for information on
causal symmetric spaces, and recall only the
necessary facts. We can always choose a maximal
abelian subspace $\fa\subset \fq_p$ which
is maximal abelian in $\fq$ and $\fp$.
Let $\Omega_A:=\Omega\cap \fa$. Then
$\Omega = \Ad (H)(\Omega_A)$, and there exists
a unique element $X_0\in \fq_p^{H\cap K}\cap \Omega_A$ such
that $\ad (X_0)$ has eigenvalues $0$, $1$ and $-1$.     
Furthermore
$\G (\Omega):=\exp (\Omega )H=H\exp (\Omega)$ is
an open semigroup, diffeomorphic to $\Omega\times H$.
 The set $\Sigma$ of roots
decomposes into   two disjoint sets:
$\Sigma_0:=\{\a\in\Sigma\mid \a (X_0)=0\}$ and
$\Sigma_n   :=\{\a\in\Sigma\mid \a (X_0)\not=0\}$.
Furthermore
$\Sigma_0$ is the set of roots of $\fa$ in $\fg^{\theta\tau}=
\fh_k\oplus \fq_p$.
 Let $W_0$ be
be the corresponding Weyl group, i.e., $W_0$ is
generated by the reflections $r_\a$, $\a\in \Sigma_0$
and --~because of our assumption $G\subset G_\C$~--
we have $W_0=N_{K\cap H}(A)/Z_{K\cap H}(A)$.
Let $\Sigma_n^+:=\{\a\in\Sigma_n\mid \a (X_0)=1\}$.
Then, since $X_0$ is $K\cap H$-invariant, it follows that
$\Sigma_n^+$ is $W_0$-invariant.    We choose our
positive set of roots such that
$\Sigma^+=\Sigma^+_0\cup \Sigma_n^+$, where $\Sigma_0^+$ is
a positive system in $\Sigma_0$. It is important to remark that
$\Sigma$ is reduced, that is 
$2\a \notin \Sigma$ for all $\a \in\Sigma$.
We denote by $\Pi_0\subset\Pi$ the set of simple
roots in $\Sigma_0^+$. Then $\Pi\setminus\Pi_0$ contains
one element.
Let $\fn_0:=\bigoplus_{\alpha\in\Sigma_0^+}\fg_\a$
and $\fn_+:=\bigoplus_{\a\in \Sigma_n^+}\fg_\a$. Similarly
we set $N_0:=\exp (\fn_0)$ and $N_+:= \exp (\fn_+)$.
Let $A_0:=\exp (\Omega_A)=\Gamma (\Omega)\cap A$. Then we have the
following theorem.

\begin{Thm}
Let the notation be as above. Then the following holds.
\begin{enumerate}
\thmlist
\item $A_0$ is $W_0$-invariant;
\item $A_0=(\overline{\bigcup_{w\in W_0}w(A^+)})^o$;
\item $\Gamma (\Omega)\subset HAN$;
\item $\Gamma(\Omega) =\Ad (H)(A_0)$;
\item $M\subset H\cap K$ and
$$H\times A\times N\ni  (h,a,n)\mapsto han\in  G$$
is a diffeomorphism onto on open subset of $G$;
\item The group $N_+$ is abelian, $G^{\tau\theta}$ normalizes $N_+$,
and $N_0\subset G^{\theta\tau}$;
\item The group $N$ can be represented as semidirect
product    $N=N_0N_+=N_+N_0$.
\end{enumerate}
\end{Thm}

Notice that (3) gives a generalized Iwasawa decomposition
for the $H$-invariant domain $\Gamma (\Omega)$
and (4) is a generalized Cartan decomposition
for  $\Gamma (\Omega)$. In particular
we note that $H$-biinvariant functions on $\Gamma(\Omega)$ are
uniquely determined by their $W_0$-invariant restriction to $A_0$.
By (5) we can define an analytic map $a_H : HAN \to A$ by
$x\in Ha_H(x)N$. For $\lambda \in \fa_\C^*$ let
\begin{equation}\label{eq-spfctGH}
\varphi_{\scriptscriptstyle{\Pi_0}}(\lambda,s)=\int_H a_H(sh)^{\lambda-\rho}\, dh \,, \quad
s\in \Gamma(\Omega),
\end{equation}
whenever $H\ni h\mapsto a_H(sh)^{\lambda -\rho}\in\C$ is
integrable. Denote by $\mathbb{D}_{G/H}(A_0)$ the
space of differential operators on $A_0$ gotten
by taking $H$-radial component of
$G$-invariant differential operators on $G/H$.
Let $\bar{N}=\theta (N)=\exp (\oplus_{\alpha\in\Sigma^+}\fg_{-\alpha})$.
Then $\bar{N}=\bar{N}_+\bar{N}_0$, with the obvious notation.
The main results in \cite{FHO} can now be stated
as follows.

\begin{Thm} Let the notation be as above.   Then the following
holds.
\begin{enumerate}
\item There exist an open convex set 
$\emptyset\not=\mathcal{E}\subset \fa_\C^*$ such
that $\varphiNCC$ is well defined for all $\l \in \mathcal{E}$ and
$a \in A_0$.
\item Assume that $\lambda \in \mathcal{E}$, then
$\varphi_{\scriptscriptstyle{\Pi_0}}(\lambda,\cdot ) |_{A_0}$ is a joint eigenfunction
for $\mathbb{D}_{G/H}(A_0)$.
\item  Let $s,t\in \Gamma (\Omega)$. Then
$H\ni h\mapsto \varphi_{\scriptscriptstyle{\Pi_0}}(\lambda,sht)\in\C$ is integrable, and
$$\int_H \varphi_{\scriptscriptstyle{\Pi_0}}(\lambda,sht)\, dh
=\varphi_{\scriptscriptstyle{\Pi_0}}(\lambda,s)
\varphi_{\scriptscriptstyle{\Pi_0}}(\lambda,t)\, .$$
\item There exists an open convex set 
$\emptyset\not=\mathcal{E}_0\subset \mathcal{E}$
such that if $\lambda \in \mathcal{E}_0$, then the integral
$$\cNCC(\lambda )
:=\int_{\bar{N}\cap HAN} a_H(\bar{n})^{-\lambda -\rho}\, d\bar{n}$$
exists and
$$\lim_{A^+\ni a\to \infty} a^{\rho -\lambda} \varphiNCC =\cNCC(\lambda)\, .$$
\end{enumerate}
\end{Thm}

Because of the decomposition $N=N_0N_+$ it follows easily that
the $c$-function can be written as a product, $\cNCC (\lambda )
=c_{\scriptscriptstyle{\Pi_0} }^+(\lambda )
c_{\scriptscriptstyle{\Pi_0}}^-(\lambda)$, where
$$c_{\scriptscriptstyle{\Pi_0}}^+(\lambda)
=\int_{N_0}a_H(\bar{n})^{-\lambda-\rho_0}\, d\bar{n}\, ,
\quad\text{with} \quad \rho_0:=\frac{1}{2}\sum_{\a\in\Sigma_0^+}\ma \a\,,$$
is the Harish-Chandra $c$-function for the Riemannian
symmetric space $G^{\theta\tau}/H\cap K$ --~even if $G^{\theta\tau}$ is not 
semisimple, Harish-Chandra's theory can still be worked out on 
$G^{\theta\tau}/H\cap K$. See the remark at the beginning of Section 
\ref{section:Rssp}~--.
Moreover, $c_{\scriptscriptstyle{\Pi_0}}^-(\lambda)$ 
is related to the geometry of $G/H$ by
$$c_{\scriptscriptstyle{\Pi_0}}^-(\lambda )=\int_{\bar{N}_+\cap HAN }
a_H(\bar{n})^{-\lambda-\rho}d\bar{n}\, .$$
Furthermore $c_{\scriptscriptstyle{\Pi_0}}^- (w\lambda )=
c_{\scriptscriptstyle{\Pi_0}}^- (\lambda )$ for
all $w\in W_0$. We remark, that we are not using the standard
notation from \cite{FHO},\cite{KO02},\cite{KO03},\cite{O97}. 
Instead we use here the
same notation as for the general case in the next section, but we
do not include the multiplicity function $m$ at
this point. The connection
to the notation in above references  is:
$c_{\scriptscriptstyle{\Pi_0}}^+ (\lambda )=c_0(\lambda)$ and
$c_{\scriptscriptstyle{\Pi_0}}^- (\lambda )=c_\Omega (\lambda )$.

Notice that the radial part of the Casimir operator is
elliptic on $A^+$ and hence
$\varphi_{\scriptscriptstyle{\Pi_0}}(\lambda,\cdot ) |_{A^+}$ is real analytic.
By Theorem \ref{th34} there exists, for generic
spectral parameters, complex numbers $c_{w,\lambda}$, $w\in W$, such
that
$\varphiNCC  =\sum_{w\in W} c_{w,\lambda} \Phi_{w\lambda}(m;a)$ for
all $a\in A^+\,.$
It was proved in \cite{O97} that the coefficients $c_{w,\lambda}$ are
given by the $c$-function, in fact
\begin{equation}\label{eq-expformspf}
\varphiNCC  =\sum_{w\in W_0} \cNCC(w\lambda)\Phi_{w\lambda} (m;a)
=c_{\scriptscriptstyle{\Pi_0}}^- (\lambda )\sum_{w\in W_0}
c_{\scriptscriptstyle{\Pi_0}}^+ (w\lambda) \Phi_{w\lambda} (m;a)
,\quad a\in A^+\, .
\end{equation}
From the integral formula for the $c$-function it follows easily that
$\cNCC(\lambda)$ extends to a meromorphic function on $\fa_\C^*$. Thus
(\ref{eq-expformspf}) implies, in particular, that $\varphiNCC $, as a function
of the spectral parameter $\lambda$, extends to a meromorphic
function on $\fa_\C^*$. Using the methods of Heckmann and Opdam it was
also proved in \cite{O97} that $\varphiNCC  $ is
smooth on all of $A_0$.
\begin{Ex}[The rank one case, \cite{FHO}] 
\label{ex:rankoneNCC}
The noncompactly causal symmetric
spaces of rank one are all locally isomorphic to $\mathrm{SO}_0(1,n)/\mathrm{SO}_0(1,n-1)$.
In this case $\Sigma=\Sigma_n=\{\alpha,-\alpha\}$. In particular
$W_0=\{1\}$. We identify
$\fa_\C^*$ with $\C$ by setting $\alpha\equiv 1$. Then
in particular $\rho=(n-1)/2$. Set $a_t=\exp (tX_0)$. Then
$\R^+\ni t\mapsto a_t\in A_0$ is a diffeomorphism. Notice, that
$\mathbb{D}_{G/H}(A_0)$ is generated by the radial part of the
Casimir element, and the corresponding differential equations
on $A_0$ are the same as for the Riemannian case. A simple
calculation (cf. \cite{FHO}) shows that
$$c_{\scriptscriptstyle{\Pi_0}}^- (\lambda ) = \kappa \frac{\Gamma (\rho)\Gamma (-\lambda
-\rho +1)}{\Gamma (-\lambda +1)}=\kappa
B(\rho, \lambda - \rho+1)$$
where $\kappa$ is a nonzero constant.
Finally we get:
$$\varphi_{\scriptscriptstyle{\Pi_0}}(\lambda,a_t)=
c_{\scriptscriptstyle{\Pi_0}}^- (\lambda )(2\cosh t)^{\lambda- \rho}
{}_2F_1\left(\frac{-\lambda-\rho}{2},
\frac{-\lambda+\rho+1}{2};1-\lambda;\frac{1}{\cosh^2t}\right)\, .$$
\end{Ex}
\begin{Ex}[The complex case, \cite{FHO}]
\label{ex:complexNCC}
The space $G_\C/G$ is NCC
if and only if $G/K$ is a bounded symmetric domain. In this case
$$\varphiNCC =
\kappa_1 c_{\scriptscriptstyle{\Pi_0}}^- (\lambda ) \, \frac{\sum_{w\in W_0}(\det w)\, a^{w\lambda}}%
{\prod_{\a\in\Sigma_0^+} \lambda_\a \; \prod_{\a\in\Sigma^+}\sinh (\a (\log a))}$$
and
$$c_{\scriptscriptstyle{\Pi_0}}^- (\lambda )=\kappa_2 \prod_{\a\in \Sigma_n^+} \lambda_\alpha^{-1}$$
where $\kappa_1$ and $\kappa_2$ are nonzero constants.
\end{Ex}

The only unknown part in (\ref{eq-expformspf}) is so far the $c_{\scriptscriptstyle{\Pi_0}}^-$-function. It was
worked out for the rank one case in \cite{FHO}, for the Cayley type spaces by
Faraut in \cite{F95}, and for $\mathrm{SL}(n,\R)/\mathrm{SO}(p,q)$ by
Graczyk \cite{G97}. The general case was then solved in \cite{KO02},
Theorem III.5, see also
\cite{KO03}.

\begin{Thm}[Product formula for the c-function]
\label{thm:comega}
For $\lambda\in\fa_\C^*$ and $\alpha\in\Sigma$ let
$\lambda_\a$ be as in (\ref{eq:la}). Then the following 
product formula holds for the function $c_{\scriptscriptstyle{\Pi_0}}^-$:
$$c_{\scriptscriptstyle{\Pi_0}}^- (\lambda )=\kappa \prod_{\a\in\Sigma_n^+} B\left(\ma/2,
-\lambda_\a-\ma/2+1\right)=
\kappa \prod_{\a\in\Sigma_n^+} \frac{\Gamma \left(\ma/2\right)%
\Gamma\left(-\lambda_\a -\ma /2 +1\right)}
{\Gamma \left(-\lambda_\a+1\right)}\,,
$$
where $\kappa$ is a nonzero constant.
\end{Thm}
Notice that, if $\ma$ is even for all $\a$, then $1/c_{\scriptscriptstyle{\Pi_0}}^- (\lambda )$ is a polynomial:
$$\frac{1}{c_{\scriptscriptstyle{\Pi_0}}^-(\lambda  )}
=\gamma \prod_{\a\in\Sigma_n^+}\prod_{j=0}^{(\ma /2)-1}(\lambda_\alpha +j)
$$
where $\gamma :=(-1)^{\sum_{\a}\ma/2} \big[\kappa \prod_{\alpha\in \Sigma_n^+} 
\Gamma (\ma /2 ))\big]^{-1}$.

Define now the \textit{spherical Fourier-Laplace transform} $(\mathcal{F}f)(\lambda )$ of
a compactly supported function on $A_0$ by
\begin{equation}\label{eq-defftr}
(\mathcal{F}f)(\lambda )=
\int_{A_0}f(a)\varphiNCC\, \delta (a)\, da\, .
\end{equation}
Then we have the following inversion formula \cite{O97}.

\begin{Thm} Let $c(\lambda )$ be the Harish-Chandra $c$-function
for the Riemannian symmetric space $G/K$, and let $\varphi_\lambda$ 
the spherical function on $G/K$. Then there exists a constant $\eta$ such that
for all $f\in C_c^\infty (A_0)$ the spherical Fourier-Laplace transform
$\mathcal{F}f$ of $f$ is inverted according to
$$
f(a)= \eta \int_{i\fa^*} (\mathcal{F}f)(\lambda )\varphi_{-\lambda }(a)
\frac{d\lambda }{\cNCC(\lambda )c (-\lambda )}\, .
$$
\end{Thm}

Paley-Wiener type theorems for the spherical Laplace-Fourier transform 
were first
considered for rank one spaces in \cite{AO01}, and have then been proven
for noncompactly causal symmetric spaces with even multiplicity  \cite{AOS} 
and for some Cayley-type spaces \cite{AU}. We do not state the
exact theorems here. In the even multiplicity case they are 
contained in the more general results for the $\Theta$-spherical 
transform, which shall be discussed in the next section. 
We would just like to point out
that the general Paley-Wiener type theorem has not been worked out so far.

\section{The $\Theta$-spherical functions}
\label{section:theta}

\noindent
The theory of $\Theta$-functions extends Heckman-Opdam's theory of
hypergeometric functions associated with root systems 
so that it encloses
special geometric instances which we have discussed: 
the Harish-Chandra's theory of spherical
functions on Riemannian symmetric spaces of noncompact type, and the
theory of spherical functions on NCC symmetric spaces.
The
$\Theta$-spherical functions,
as the Heckman-Opdam's hypergeometric functions,
are special functions associated with root
systems. Because of Examples \ref{ex:rankone} and \ref{ex:rankoneNCC},
they can be considered as geometrically motivated multivariable
generalizations of the Jacobi functions.

Both spherical functions on Riemannian and NCC symmetric spaces
admit expansion with respects to the Harish-Chandra series. 
This suggests, exactly as in the case of hypergeometric functions
associated with root systems, to use linear combinations of
Harish-Chandra series for defining their common generalizations.
Comparing the expansion (\ref{eq-expformspf}) to 
(\ref{eq:orHCexpansion}), one notices
that the summation over the Weyl group $W$ is replaced by a summation
over its parabolic sugroup $W_0$. The coefficients
for the Harish-Chandra series are also modified. 
A unified theory can therefore be constructed by considering suitable
linear combinations of Harish-Chandra series over parabolic
subgroups of $W$.

For semplicity of exposition we shall assume in the following that the
root system $\Sigma$ is reduced, i.e. that $2\a \notin \Sigma$ for all
$\a \in \Sigma$. This condition is for instance always satisfied by the
root systems of NCC symmetric spaces. We refer to \cite{OP1}, \cite{P1}
and \cite{P2} for the theory of $\Theta$-spherical functions on 
general root systems and for further details. We keep the notation introduced in
Section \ref{section:HO}.
 
Let $\Theta$ be an arbitrary subset of positive simple roots
in a fundamental system $\Pi \subset \Sigma^+$, and $\pedtheta{W}$
the parabolic subgroup of $W$ generated by the reflections $r_\a$,
$\a \in \Theta$. Then $\pedtheta{W}$  is the Weyl group of the root system
$\rootstheta:=\Z \Theta \cap \Sigma$. 
Set $\rootstheta^+:=\rootstheta \cap \Sigma^+$. Notice that 
$\pedtheta{W}(\Sigma^+\setminus \rootstheta^+)
\subset \Sigma^+\setminus \rootstheta^+$.
We refer to \cite{Hu} for additional information on parabolic subgroups of
Weyl groups.

We introduce the $c$-functions
\begin{align*}
\cthetampl&:=\prod_{\a \in \rootstheta^+} \frac{\Gamma(\la)}
                   {\Gamma(\la+\ma/2)}\,,\\
\cthetamml&:=\prod_{\a \in \Sigma^+ \setminus \rootstheta^+}
               \frac{\Gamma(-\la-\ma/2+1)}
                {\Gamma(-\la+1)}\,,
\end{align*}
with the conventions that products over empty sets are equal to $1$, i.e.
$$c_\emptyset^+=\cPi^{-} :=1.$$
Notice that (up to a constant depending on $m$ but not on $\l$)
the function $\pedtheta{c}^+$ is Harish-Chandra's $c$-function for the root
system $\rootstheta$, whereas $\pedtheta{c}^-$ is modelled on the function 
$c^-_{\scriptscriptstyle{\Pi_0} }$
of Theorem \ref{thm:comega}. Observe also that
$\pedtheta{c}^-(m;\l)$ is a 
$\pedtheta{W}$-invariant function of $\l \in \frakacs$. 

When all multiplicities $m_\a$ are even, 
all these $c$-functions are reciprocals of polynomials.
Indeed, in this case
\begin{align*}
\cthetampl&:=\prod_{\a \in \rootstheta^+} \prod_{h=0}^{\ma/2-1} 
\frac{1}{\la+h}\,,  \\
\cthetamml&:=(-1)^{d(\Theta,m)}  \prod_{\a \in \Sigma^+\setminus \rootstheta^+} \;
 \prod_{h=0}^{m_\a/2-1} \frac{1}{\la+h}\,,
\end{align*}
where
\begin{equation}
  \label{eq:dmtheta}
d(\Theta,m):=\frac{1}{2} \sum_{\a\in \Sigma^+\setminus \rootstheta^+}\ma\,.
\end{equation}

\begin{Def} \label{def:thetasph}
Let $\Theta\subset \Pi$ and
let $U^+$ be a tubular neighborhood of $A^+$ in $\complex{A}$
on which, for generic $\l\in \frakacs$, 
each series defining $\Phi_{w\l}(m;a)$ converges to a holomorphic function of $a$
(cf. comments on p. \pageref{Phiholo}).
The function on $U^+$ defined for generic $\l \in \frakacs$ by
\begin{equation}
  \label{eq:thetaphi}
  \pedtheta{\varphi}(m;\l,a):=\pedtheta{c}^-(m;\l) \; \sum_{w \in \pedtheta{W}}
\pedtheta{c}^+(m;w\l) \Phi_{w\l}(m;a), \qquad a \in U^+
\end{equation} 
is called the \emph{$\Theta$-spherical function of spectral parameter $\l$}.
\end{Def}

As a linear combination of the Harish-Chandra series $\Phi_{w\l}(m;a)$, the 
$\Theta$-spherical function of spectral parameter $\l$ is, by construction,
 a solution of the
hypergeometric system (\ref{eq:hypereq}).

We now illustrate how to recover our motivating examples from the context of $\Theta$-spherical
functions.

\begin{Ex}
\begin{enumerate}
\thmlist
\item 
If $m=0$, then $\Phi_\l(0;a)=a^\l$. We thus obtain the Euclidean case (with different
symmetries induced by the choice of $\Theta$).
\item
When $\Theta=\Pi$, then $\pedtheta{\varphi}(m;\l,a)$ coincides --~up to the normalization factor
$\pedPi{c}^+(m;\rho(m))$~-- with Heckman-Opdam's hypergeometric function of
spectral parameter $\l$, hence with Harish-Chandra's spherical function of
spectral parameter $\l$ when the triple $(\frak a,\Sigma,m)$ is geometric.
\item
Suppose $(\frak a,\Sigma,m)$ corresponds 
to the Riemannian dual of an NCC symmetric space $G/H$
according to Section \ref{se-nonRsp}. 
If $\Theta=\Pi_0$ is the set of positive compact simple
roots, then $\ped{\Pi_0}{\varphi}(m;\l,a)$  coincides --~up to the normalization factor
$\ped{\Pi_0}{c}^+(m;\rho(m))$~-- with the spherical function 
$\varphi_{\scriptscriptstyle{\Pi_0}}(\l,a)$
on $G/H$. Notice also that $\ped{\Pi_0}{W}=W_0$.
\end{enumerate}  
\end{Ex}

As in the theory of Heckman and Opdam, the first crucial problems are to understand the domain
on which the $\Theta$-spherical functions extend as functions of both $\l$ and $a$, and
to study the nature of this extension. Notice that in the NCC case of Section \ref{se-nonRsp}
it was proved that the spherical functions are meromorphic in $\l \in \frakacs$, but there was
no description of the $\l$-singular set. The $\l$-singular set of the
$\Theta$-spherical functions 
--~and hence also of the spherical functions on NCC symmetric
spaces~-- turns out to be described by the numerator $\pedtheta{n}^-$
of the function $\pedtheta{c}^-$:
\begin{equation*}
\pedtheta{n}^-(m;\l):=
\prod_{\a \in \Sigma^{+}\setminus \rootstheta^+} 
\Gamma\left(-\la-\frac{\ma}{2}+1\right).
\end{equation*} 
In the $a$-variable the situation is very similar to the one depicted for the special 
case of spherical functions on NCC symmetric spaces.
The domain of extension will be
\begin{equation}
  \label{eq:eq:atheta}
  \pedtheta{A}:=\exp(\pedtheta{\frak a}) \quad\text{with} \quad 
\pedtheta{\frak a}:=\Big(\overline{\bigcup_{w\in\pedtheta{W}} w(\frak a^+)}
\Big)^o\,.
\end{equation}
The proof of the following theorem can be found in \cite{OP1}, Theorem 8.3,
or \cite{P1}, Theorem 3.5.

\begin{Thm} 
  \label{thm:entire}
Let the notation be as above. 
Then there exists a $\pedtheta{W}$-invariant 
tubular neighborhood $\pedtheta{U}$
of $\pedtheta{A}$ in $\complex{A}$ such that the function 
\begin{equation*}
(\l,a) \longmapsto \frac{\pedtheta{\varphi}(m;\l,a)}
{\pedtheta{n}^-(m;\l)}
\end{equation*}
extends as a $\pedtheta{W}$-invariant
holomorphic function of $(\l,a) \in \frakacs \times \pedtheta{U}$.
\end{Thm}

Observe that (if $\Sigma$ is reduced as we are assuming here)
the $\l$-singularities are  at most of first order and 
located along a locally finite (in general
infinite) family of complex affine hyperplanes in $\frakacs$.
We shall see shortly that the situation greatly simplifies 
when all multiplicities are even.

\begin{Def} \label{def:fourier}
Let $\Theta \subset \Pi$.
The \emph{$\Theta$-spherical
transform} of $f \in C^\infty_c(\pedtheta{A})^{\pedtheta{W}}$ 
associated with the data $(\frak a, \Sigma, m)$ is the $\pedtheta{W}$-invariant
function $\pedtheta{\mathcal F}f(m)$ on $\frakacs$ defined for $ \l \in \frakacs$ by
\begin{equation}
\pedtheta{\mathcal F}f(m;\l):=
\frac{1}{\abs{\pedtheta{W}}} \, \int_{\pedtheta{A}}
f(a) \, \pedtheta{\varphi}(m;\l,a) \delta(m;a)\; da 
=\int_{A^+} f(a) \pedtheta{\varphi}(m;\l,a) \delta(m;a)\; da\,,
\end{equation}
where $da$ is a suitable normalization of the Haar measure on $A$.
\end{Def}

For every $\l \in \frakacs$ let $\pedtheta{E}(m;\l)$ be the 
$\pedtheta{W}$-invariant function 
on $\pedtheta{A}$ defined by requiring that the equality
\begin{equation} \label{eq:ETheta}
\pedtheta{E}(m;\l,a)=
\frac{\pedtheta{c}^-(m;\l)\pedtheta{c}^+(m;\l)}{\pedPi{c}^{+}(m;\l)} \;
\pedPi{\varphi}(m;\l,a), \qquad a \in A^+\,.
\end{equation}
Then $\pedtheta{E}(m;\l,a)$ is a $\pedtheta{W}$-invariant 
meromorphic function of $\l\in \frakacs$.
Notice that, if  $m$ is an even multiplicity function, then
$\pedtheta{E}(m;\l,a)=\pm \pedPi{\varphi}(m;\l,a)$ for all $\Theta$.
The inversion of the $\Theta$-spherical transform on 
$C^\infty_c(\pedtheta{A})^{\pedtheta{W}}$ is provided by the following theorem.
See \cite{P2}, Theorem 4.5.

\begin{Thm} \label{thm:inv}
Let $f \in C^\infty_c(\pedtheta{A})^{\pedtheta{W}}$. Then for all
$a \in \pedtheta{A}$ we have
\begin{equation} \label{eq:inversion}
f(a)= \kappa\; \frac{\abs{W}}{\abs{\pedtheta{W}}}
\int_{i\frak a^*} \pedtheta{\mathcal F}f (m;\l) \, \pedtheta{E}(m;-\l,a)
 \abs{\pedtheta{c}(m;\l)}^{-2}\;  d\l\,,
\end{equation}
where $d\l$ is a suitable normalization of the Lebesgue
measure on $i\frak a^*$ and
$\kappa$ is a positive constant depending on
the normalizations of the measures.
\end{Thm}

A Paley-Wiener theorem for the $\Theta$-spherical transform is a
characterization of the image under $\pedtheta{\mathcal F}(m)$ of the set 
$C^\infty_c(\pedtheta{A})^{\pedtheta{W}}$. 
A particularly interesting 
higher-rank situation in which this theorem has been worked out is
the even multipicity case.

The \textit{even multiplicity case} corresponds to triples $(\frak a,\Sigma,m)$ 
for which $\Sigma$ is a reduced root system and $m_\a \in 2\N_0$ for all
$\a \in \Sigma$.
In the geometric case of Riemannian symmetric spaces of noncompact type,
this situation corresponds to spaces $G/K$ 
in which the Lie algebra $\frak g$ of $G$ possesses a unique conjugacy
class of Cartan subalgebras. Moreover, as already observed in 
Examples \ref{ex:complex} and \ref{ex:complexNCC},  
multiplicity functions with constant value $2$ 
correspond to  $\frak g$'s admitting a complex structure.

The class of symmetric spaces for which the theory of $\Theta$-spherical
functions is relevant consists of the $K_\e$ symmetric spaces of Oshima 
and Sekiguchi \cite{OS}. The most important examples among them
are our basic examples, that is the  Riemannian symmetric spaces of noncompact 
type and the NCC symmetric spaces.
Let $G/H$ be a $K_\e$ symmetric space with even multiplicity functions
with $G$ a connected, noncompact, simple Lie group. 
Then, according to classification (\cite{OS}; see also 
Appendix \ref{app:classification} below),
either $\frak g$ admits complex structure or the root system of 
$G/K$ is of type $A_n$.
In the latter case, $G/H$ is either Riemannian or NCC. Furthermore,
all irreducible $K_\e$ spaces are determined by a
signature, corresponding to the choice of a set $\Theta\subset \Pi$
of positive simple roots so that either $\Theta=\Pi$ (Riemannian case)
or $\abs{\Pi \setminus \Theta}=1$. 
We shall therefore state the Paley-Wiener theorem for the
$\Theta$-spherical transform in even multiplicites with the 
additional assumptions on
$(\frak a,\Sigma,m)$ that $\Sigma$ 
is a root system of type $A_n$
and that $\Theta$ consists of at least $n-1$ positive simple roots. 
As just observed, this is the most
significant case from the point of view of the harmonic analysis on 
symmetric spaces with even multiplicities. By only dealing with the above
situation, we can explain the proving methods 
without entering into the technical difficulties which appear in the
general even multiplicity case. 
For the latter, we refer the reader to \cite{OP3}.

One of the simplifications occurring in the case of $A_n$
is that there is only one Weyl group orbit in $\Sigma$.
This forces the multiplicity function
$m$ to be constant. In the following we shall denote by the same symbol
$m$ the constant even multiplicity function and its value. 

Before stating the Paley-Wiener therorem, let us add some important 
properties of the $\Theta$-spherical functions in the even multiplicity
case.
Compared with the general case, the even multiplicity case (not necessarily 
$A_n$) is greatly simplified by the fact that the $\l$-singularities are
located only along a {\em finite} family of complex affine hyperplanes. 
Moreover -- and this will turn out to be one of the key ingredients in proving
the Paley-Wiener theorem -- there are explicit formulas
relating, by means
of differential operators,
the $\Theta$-spherical functions to exponential functions. 
The main tool for this are 
\textit{Opdam's shift operators}.

Let $\C[\complex{A}]$ be the algebra of regular functions on $\complex{A}$ 
introduced in Section \ref{section:HO}, and let $\Delta$ be the Weyl
denominator as in (\ref{eq:Weylden}). The algebra 
$\C_\Delta[\complex{A}]=\bigcup_{k\in\Z} \, \Delta^k  \C[\complex{A}]$
is the localization of $\C[\complex{A}]$ along $\Delta$. Shift 
operators are $W$-invariant differential operators on $A$ with coefficients in
$\C_\Delta[\complex{A}]$. They have a characteristic asymptotic expansion 
on $A^+$ and commute with the operator $L(m)$.
For every even multiplicity function $m$ there are shift 
operators ``lowering'' the
Harish-Chandra series $\Phi_\l(m;a)$ to the
exponential function $a^\l$ (which is the Harish-Chandra series for $m=0$),
or, vice versa, ``raising'' $a^\l$ to $\Phi_\l(m;a)$.

\begin{Lemma}[Opdam] \label{lemma:shift}
  Let $m \in 2\N_0$ be even. Then there exists unique shift operators
$D_-(m)$ and $D_+(m)$ so that for all $(\l,a) \in 
(\frakacs\setminus P) \times A^+$
\begin{align}
D_-(m)\Phi_\l(m;a)&=\frac{1}{\pedPi{c}^+(m;\l)} \; a^\l\,,
\label{eq:DminusHC}\\
D_+(m)a^{\l}&= \frac{1}{\pedPi{c}^+(m;-\l)} \; \Phi_\l(m;a).
\label{eq:DplusHC}
\end{align}
\end{Lemma}

The shift operators $D_{\pm}(m)$ of Lemma \ref{lemma:shift} can be 
constructed by composing fundamental shift operators
of shifts $\pm 2$. Let  $G_+(m)$ and $G_-(m)$ respectively 
denote the fundamental shift operator  
mapping $\Phi_\l(m;a)$ to $\Phi_{\l}(m+2;a)$ and the
fundamental shift operator  
mapping $\Phi_\l(m;a)$ to $\Phi_{\l}(m-2;a)$. Then $G_+(m)$
is related to the formal adjoint of   
$G_-(m)$ by the relation
\begin{equation}
\label{eq:gpluseminus}
G_+(m):=\Delta^{-(2+m)}\circ G^*_-(m+2)\circ \Delta^m\,.
\end{equation}
We refer to Proposition 4.4 in \cite{HeckBou}
for general explicit formulas.
In the geometric case of Riemannian symmetric spaces, 
$D_+(m)$ is a constant multiple of the
differential operator inverting the Abel transform. Further results 
and explicit formulas have been obtained in this context by Beerends
(\cite{BeeComp}, \cite{BeeCan}), Hba \cite{Hba}, and Vretrare \cite{Vre}.

\begin{Ex}[The root system $A_1$]
  \label{ex:shiftA1}
In the rank-one $A_1$ case there is only one positive root $\a$. With the
identifications introduced in Example \ref{ex:rankone}, 
the shift relations of Lemma \ref{lemma:shift}
correspond to the classical differentiation formula for the
hypergeometric function
\begin{equation*}
\frac{d}{dz} \hyper{a}{b}{c}{z}=\frac{ab}{c} \, \hyper{a+1}{b+1}{c+1}{z}
\end{equation*}
(cf. e.g.  \cite{Er}, 2.8(20)).
Seen as differential operators on $\liecomplex{a}\equiv \C$ 
the elementary shift operators of shift 2 are therefore
\begin{align*}
G_-(m)&=\Delta(z)\frac{d}{dz}+(m-1)(e^z+e^{-z}), \\
G_+(m)&=-\frac{1}{\Delta(z)}\;\frac{d}{dz}  
\end{align*}
with $\Delta(z):=e^z-e^{-z}$. 
Since $G_+(m)$ is independent of $m$,  it follows in particular that
$$D_+(m)=G_+(2)^{m/2}.$$ 
\end{Ex}

\begin{Ex}[The root system $A_2$]
  \label{ex:shiftA2} 
The fundamental shift operators $G_+(m)$
for the root system $A_2$ were determined by Beerends \cite{BeeComp}. 
They are
\begin{equation*}
G_+(m)=\Delta^{-1} \Big[ \prod_{\a\in \Sigma^+} \partial_\a +
  \frac{m-2}{2} \sum_{\a \in \Sigma^+} \eta(\a) \circ \partial_\a
  \circ \coth \a \circ \partial_\a \Big]\,
\end{equation*}
where
\begin{equation*}
\eta(\a)=\prod_{\b\in \Sigma^+ \setminus\{\a\}} \inner{\a}{\b}.
\end{equation*}
The shift operators $G_-(m)$ can then be deduced from
(\ref{eq:gpluseminus}). 
Notice that the formula semplifies in the complex case, since
\begin{equation*}
G_+(0)=\Delta^{-1} \; \prod_{\a\in \Sigma^+} \partial_\a.
\end{equation*}
Composition yields for instance 
\begin{multline*}
D_+(4)=G_+(2)\circ G_+(0)
 =\Delta^{-2}\left[\partial_{\a}\partial_{\b}\partial_{\gamma} -{\sum}' 
  \abs{\a}^2 \coth(\a)\partial_{\beta}
 \partial_{\gamma} + {\sum}' \abs{\a}^2 \abs{\b}^2 \coth(\a) \coth(\b)
\partial_{\gamma} \right. \\
 \left. -\abs{\a}^2\abs{\b}^2 \abs{\gamma}^2 \coth(\a)\coth(\beta)\coth(\gamma)
 -\frac{\abs{\a}^2\abs{\beta}^2 \abs{\gamma}^2}{2\Delta}\right] \circ
  \partial_{\a}\partial_{\b}\partial_{\gamma}\,, 
\end{multline*}
where ${\sum}'$ denotes the sum over the cyclic permutations of the
three positive roots $\a,\b,\gamma=\a+\b$.
The operator $D_+(4)$ was first computed by Hba \cite{Hba} 
as the inverse of the 
Abel transform on the Riemannian symmetric space
$\SL(3,\quat)/{\rm Sp}(3)$.
\end{Ex}

Lemma \ref{lemma:shift} together with Definition \ref{def:thetasph}
immediately yields explicit formulas for the $\Theta$-spherical functions, 
at least on $A^+$. The difficulty is to extend these formulas to the
domain $\pedtheta{A}$, since the shift operators have singular 
coefficients on the set $\{\Delta=0\}$. The definition of
shift operators ensures that these singularities 
can be cancelled by multiplication by a suitable power
$\Delta^k$, but it gives no estimates on $k$. Direct estimates
from Heckman's explicit formulas for the fundamental shift operators are
hard. First of all, these formulas are rather involved in 
high rank. Moreover, in the composition of shift
operators of step 2 many simplifications occur. 
In fact the exponent $k$ is lower than one could guess
from the form of the
single step-2 parts. 
This occurs for instance for the operator $D_+(4)$
in the above Example \ref{ex:shiftA2}.

An argument for estimating the value of $k$ has been presented
in \cite{OP3}, Theorem 4.10.
It does not rely on explicit formulas for the shift operators, but on general
estimates for the Harish-Chandra series and their derivatives.
It states that $k$ is smaller or equal to $m$. This in particular
implies that $D_m:=\Delta^m D_+(m)$ extends as $W$-invariant
differential operator on 
$A$ with coefficients in $\C[\complex{A}]$.

The regularity properties of the $\Theta$-spherical functions in
even multiplicites are collected in the following theorem. We state it
for the case in which $m \in 2\N_0$ is a constant multiplicity function,
but it can be extended with minor modifications to arbitrary even
multiplicities. We refer the reader to \cite{OP2} and \cite{OP4} for
more information.

\begin{Thm} \label{thm:thetaspheven}
Suppose $m\in 2\N_0$. Let $\Theta \subset \Pi$ and define
\begin{equation}
\label{eq:e}
\pedtheta{e}^-(m;\l):=
 \prod_{\a \in \Sigma^+\setminus \langle\Theta\rangle^+}
\; \prod_{k=-m/2+1}^{m/2-1} (\la-k)
\end{equation}
for $\l \in \frakacs$. We convene that empty products are equal to one.
\begin{enumerate}
\thmlist
\item  \label{item:poles}
There is a $\pedtheta{W}$-invariant tubular 
neighborhood $\pedtheta{U}$ in $\complex{A}$
of $\pedtheta{A}$ such that the function
\begin{equation*}
\pedtheta{e}^-(m;\l)\;\pedtheta{\varphi}(m;\l,a)
\end{equation*}
extends as a $\pedtheta{W}$-invariant holomorphic function of
$(\l,a) \in \frakacs \times \pedtheta{U}$.
\item
There is a $W$-invariant differential operator with coefficients in $\C[\complex{A}]$,
namely $D_m=\Delta^m D_+(m)$, so that for all
 $(\l,a) \in \frakacs \times \pedtheta{U}$
\begin{equation}
  \label{eq:formulabis}
  \Delta^m(a) \pedtheta{\varphi}(m;\l,a)=
(-1)^{d(\Theta,m)}
\Big[{\prod_{\a\in\Sigma^+} \prod_{k=0}^{m/2-1} (k^2-\la^2)}\Big]^{-1} \;
D_m \Big(\sum_{w \in \pedtheta{W}}a^{w\l}\Big)
\end{equation}
with $d(\Theta,m)$ as in (\ref{eq:dmtheta}).
\item
If $\Theta=\Pi$, then $\Delta^m(a) \pedPi{\varphi}(m;\l,a)$ extends as a
holomorphic function on $\frakacs \times \pedPi{U}$ by means of the formula
\begin{equation}\label{eq:formulaPi}
 \Delta^m(a) \pedPi{\varphi}(m;\l,a)=
\Big[{\prod_{\a\in\Sigma^+} \prod_{k=0}^{\ma/2-1} (k^2-\la^2)}\Big]^{-1} \;
D_m \Big(\sum_{w \in W}e^{w\l(\log a)}\Big).
\end{equation}
If, moreover, $\l \in P$, then $\Delta^m(a) 
\pedPi{\varphi}(m;\l,a)$ extends
by (\ref{eq:formulaPi}) as $W$-invariant entire function on $\complex{A}$.
\end{enumerate}
\end{Thm}

Observe that in the complex case 
\begin{equation}
D_+(m)=\Delta^{-1} \; \prod_{\a\in \Sigma^+} \partial_\a\,.
\end{equation}
In the geometric case, Formula (\ref{eq:formulaPi}) reduces therefore to the
classical formula by Harish-Chandra of Example \ref{ex:complex}, whereas
(\ref{eq:formulabis}) with $\Theta=\Pi_0$ reduces to the results
of \cite{FHO} 
explained in Example
\ref{ex:complexNCC}.

The description of the image of the compactly supported functions 
under the $\Theta$-spherical transforsm is given by means of the following
definition (which is given for the general even multiplicity case).

\begin{Def}[Paley-Wiener space] \label{def:pwspacetheta}
Let $m$ be an even multiplicity function on the root system $\Sigma$, 
and let $\Theta \subset \Pi$ be a fixed
set of positive simple roots. Let $C$ be a compact, convex and 
$\pedtheta{W}$-invariant subset of $\pedtheta{\frak a}$.
The {\em Paley-Wiener space} $\pwthetamC$ is the space of all 
$\pedtheta{W}$-invariant meromorphic functions $g:\frakacs \to \C$
satisfying the following properties:
\begin{enumerate}
\renewcommand{\labelenumi}{\theenumi.}
\item  \label{item:pwspacethetauno}
$\pedtheta{e}^-(m;\l) g(\l)$ 
is a rapidly decreasing entire function of exponential type $C$,
that is 
for every $N \in \N$
there is a constant $C_N\geq 0$ such that
$$
\abs{\pedtheta{e}^-(m;\l) g(\l)} \leq C_N (1+\abs{\l})^{-N} e^{q_C(\Re\l)}
$$
for all $\l \in \frakacs$.
\item \label{item:pwspacethetadue}
The function 
\begin{equation}\label{eq:Pav}
\Pavtheta g(\l):=\sum_{w \in \pedtheta{W}\setminus W} g(w\l)
\end{equation}
 extends to an entire function on $\frakacs$.
\end{enumerate}
\end{Def}

Condition \ref{item:pwspacethetadue} is automatically 
satisfied in the Euclidean case 
$m=0$, in the complex case $m=2$, and 
when $\Theta=\Pi$. 
Indeed, $\pedtheta{e}^-\equiv 1$ in the Euclidean case and when
$\Theta=\Pi$. For the complex case, 
observe that  Condition \ref{item:pwspacethetauno} in Definition \ref{def:pwspacetheta}
implies that $\Pavtheta g$ is $W$-invariant and has at most first order singularities 
along each hyperplane $\{\l \in \frakacs\mid \l_\a=0\}$ with 
$\a \in \Sigma^+$. 
This is only possible when 
the singularities are in fact removable.

\begin{Thm}[Paley-Wiener theorem]\label{thm:pw}
Let $\Sigma$ be a root system of type $A_n$ with multiplicity function 
$m$, and let $\Theta \subset \Pi$ 
be a set of positive simple root so that $\abs{\Theta}\geq n-1$.
Suppose $C$ is a compact, convex and $\pedtheta{W}$-invariant subset of
$\pedtheta{\frak a}$.
Then the $\Theta$-spherical transform $\pedtheta{\mathcal F}(m)$ 
maps $C^\infty_c(C)^{\pedtheta{W}}$ bijectively onto $\pwthetamC$.
\end{Thm}

Notice first that the case $\Theta=\Pi$ is contained in the Paley-Wiener
theorem for the Opdam transform (Theorem \ref{thm:pwopdam}). 
An elementary proof of the even multiplicity case with $\Theta=\Pi$ 
can also be found in \cite{OP3}, Corollary 10.2. 
In the following we shall therefore 
suppose that  $\Pi \setminus \Theta$ consists of a single element, 
say $\Pi \setminus \Theta=\{\b\}$.

The fact that the $\Theta$-spherical transform maps 
$C^\infty_c(C)^{\pedtheta{W}}$ 
into the Paley-Wiener space $\pwthetamC$ depends mainly on the 
explicit formula for the
spherical functions. 
Suppose $f \in C^\infty_c(C)^{\pedtheta{W}}$.
Observe that $\pedtheta{e}^-(m;\l)\pedtheta{\mathcal F}f(m;\l)$ is entire by 
Theorem \ref{thm:thetaspheven}. 
For $\l \in \frakacs$ set
\begin{equation}
\label{eq:eplus}
\pedtheta{e}^+(m;\l):=(-1)^{m\abs{\rootstheta^+}/2} \;
 \prod_{\a \in \Sigma^+\setminus \langle\Theta\rangle^+}
\; \prod_{k=-m/2+1}^{m/2-1} (\la-k)
\end{equation}
(with the convention that empty products are equal to $1$). 
The factor $\Delta^m$ used in the definition of $D_m$ 
appears as density of the measure inside the $\Theta$-spherical 
transform. 
Hence, Theorem \ref{thm:thetaspheven} and the $W$-invariance of $D_m$ imply 
\begin{equation} \label{eq:directPW}
\pi(\l)\pedtheta{e}^+(m;\l) \pedtheta{e}^-(m;\l) \pedtheta{\mathcal F}f(m;\l)=
\big[\mathcal F_A (D_m^*f)\big](\l).
\end{equation}
The regularity of $D_m$ ensures that $D_m^*f$ 
is a $\pedtheta{W}$-invariant smooth 
function with compact support in $\exp C$. The 
classical Paley-Wiener theorem for the Euclidean 
Fourier transform $\mathcal F_A$ (see Section \ref{section:euclPW})
implies therefore that the function in (\ref{eq:directPW}) 
belongs to $\PW(C)$. 
Now, a classical result of Malgrange states for an entire function 
$F$ and a polynomial 
$p$ that $pF \in PW(C)$ if and only if $F\in PW(C)$. 
See e.g.  \cite{He3}, Lemma 5.13. This proves the necessity of Condition
\ref{item:pwspacethetauno}  for 
$\pedtheta{\mathcal F}f$. 

Condition \ref{item:pwspacethetadue} follows for the Paley-Wiener theorem for 
the Opdam transform. 
Indeed, every function $f \in \ccitheta$ can be uniquely extended to a $W$-invariant function 
$\pedPi{f} \in C_c^\infty(A)^W$. 
The definition of $\Theta$-spherical functions yields the relation
 \begin{equation*}
  \pedPi{\varphi}(m;\l,a)=(-1)^{d(\Theta,m)} \sum_{\pedtheta{W}\setminus W} 
  \pedtheta{\varphi}(m;w\l,a),
    \end{equation*}
where $d(\Theta,m)$ is as in 
(\ref{eq:dmtheta}). 
Consequently 
\begin{equation*}
  \big(\pedPi{\mathcal F}\pedPi{f}\big)(m;\l)=
(-1)^{d(\Theta,m)} \sum_{\pedtheta{W}\setminus W} 
  \big(\pedtheta{\mathcal F}f)(m;w\l)=
(-1)^{d(\Theta,m)}  \big(\Pavtheta  \pedtheta{\mathcal F}f)(m;\l).
\end{equation*}

One should remark that in the above arguments the assumptions on 
$\Sigma$ and $\Theta$ play no role. Indeed, the 
$\Theta$-spherical transform maps $C^\infty_c(C)^{\pedtheta{W}}$ 
into $\pwthetamC$ for every root system $\Sigma$, every even multiplicity function 
on $\Sigma$, and every choice of $\Theta \subset \Pi$.    
The assumptions made in Theorem \ref{thm:pw} 
enter only in the proof of the surjectivity of the transform, which
we now outline.

The inversion formula and the explicit formulas 
for the $\Theta$-spherical functions suggest the following definition of wave packets.
The \textit{wave packet} of  $g\in \pwthetamC$ is the function 
$\mathcal I g=\mathcal I g(m):A \to \C$ defined by
\begin{equation}\label{eq:wave}
\big(\mathcal I g\big)(a)
:=\int_{i\frak a^*} g(\l) \pedPi{\varphi}(m;-\l,a) \abs{\pedPi{c}^+(m;\l)}^{-2}\; d\l 
\end{equation}
The {\em $\Theta$-wave-packet}
of $g$ is the function on $\pedtheta{A}$ obtained by restriction
of $\mathcal I g$ to $\pedtheta{A}$, that is 
\begin{equation}
  \label{eq:thetawave}
\pedtheta{\mathcal I} g = \mathcal I g \circ \pedtheta{\iota},   
\end{equation}
where $\pedtheta{\iota}:\pedtheta{A} \hookrightarrow A$ is the inclusion map.

Notice that $\mathcal I g=\frac{\abs{\pedtheta{W}}}{\abs{W}} \mathcal I \Pavtheta g$.
Hence, for some constant $\kappa$, 
\begin{equation*}
  \Delta^m(a) \mathcal I g(a)=\kappa D_m \mathcal F^{-1}_A \big(\Pavtheta g\big)(a).
\end{equation*}
Condition 2 in Definition \ref{def:pwspacetheta} and the classical Paley-Wiener theorem 
therefore ensure that the wave packet  $\mathcal I g$ has support contained in 
$\exp\big(\conv(W(C))\big)$, where $\conv(W(C))$ denotes the convex hull of the
Weyl group orbit of $C$. 

The next task is to show that $\pedtheta{\mathcal I} g$ is smooth and 
compactly supported in $\pedtheta{A}$.
On the Lie algebra level, we need to separate points 
in the interior of $\pedtheta{\frak a}$ from points on its boundary. 
This can be done by 
means of suitable elements in 
\begin{equation*}
  \pedtheta{\frak a}^*(m)=\{\l \in \pedtheta{\frak a}^* \mid 
\text{$\la\geq m/2$ for all $\a \in \Sigma^+ \setminus \rootstheta^+$}\}.
\end{equation*}
Here $\pedtheta{\frak a}^*:=\{\l \in \frak a^* \mid
\text{$\l(H)\geq 0$ for all $H \in \pedtheta{\frak a}$}\} 
   =\sum_{\a \in \Sigma^+ \setminus \rootstheta^+} \R_0^+ \a$  is the dual cone of 
$\pedtheta{\frak a}$.
The set $\pedtheta{\frak a}^*(m)$ is introduced because it is a ``large'' closed
subset of $\pedtheta{\frak a}^*$ which is ``away'' from the possible singularities of 
every $g \in \pedtheta{\PW}(m;C)$. Indeed, each $g \in \pedtheta{\PW}(m;C)$ 
is holomorphic in a neighborhood of the 
convex set $i\frak a^*-\pedtheta{\frak a}^*(m)$.
Furthermore, for every $N \in \N$, there is a 
constant $C_N>0$ such that for all 
$\l \in i\frak a^*$ and $\mu \in -\pedtheta{\frak a}^*(m)$
\begin{equation}
  \label{eq:estgaway}
  \abs{g(\l+\mu)} \leq C_N(1+\abs{\l})^{-N} e^{q_C(\mu)}. 
\end{equation}
This allows to shift the contour of integration and get 
for all $\mu \in -\pedtheta{\frak a}^*(m)$ and $a \in A$,
\begin{equation} \label{eq:partialsum}
\Delta^m(a) \mathcal I g(a)
 =\sum_{w\in W} D_m \int_{i\frak a^*} g(\l+\mu) a^{-w(\l+\mu)} \; d\l\,.
\end{equation}
The assumption that $\Sigma$ is of type $A_n$ ensures that 
$\inner{\beta}{\a} \geq 0$ for all $\a \in \Sigma^+ \setminus \rootstheta^+$.
This is used to prove that, for $a \in \pedtheta{A}$,  
only the summands of (\ref{eq:partialsum}) 
corresponding to $w \in \pedtheta{W}$ 
are nonzero. 
Hence for all $\mu \in -\pedtheta{\frak a}^*(m)$ and $a\in \pedtheta{A}$
we have
\begin{equation} \label{eq:inversiondue} 
\Delta^m(a) \pedtheta{\mathcal I} g(a)
 =\abs{\pedtheta{W}} \;  D_m \int_{i\frak a^*} 
g(\l+\mu) a^{-(\l+\mu)} \; d\l\,.
\end{equation}
Formula (\ref{eq:inversiondue}) yields the shift of contour of integration 
to show that the support of $\pedtheta{\mathcal I} g$ is a compact subset of 
$\pedtheta{A}$. 

The final step, which proves that the support of  $\pedtheta{\mathcal I} g$
is indeed contained in $\exp C$, requires a certain application of
Holmgren's theorem. Holmgren's theorem has been employed in the proof of 
Paley-Wiener type theorems also in \cite{BS}, but, 
to be able to apply it to our situation, several 
adjustments are required. 

The basic difficulty in working with $\pedtheta{\mathcal I}g$ is due to the 
possible $\l$-singularities of $g$. 
Of course one would like to replace $g$ with
$\pedtheta{e}^-(m;\l)g(\l)$, which is entire. 
For this, the trick is to use suitable differential
operators. The polynomial 
\begin{equation*}
  q(m;\l):=\prod_{\a \in \Sigma}\; \prod_{k=-m/2+1}^{m/2-1} (\l_\a-k)
\end{equation*}
is divided by $\pedtheta{e}^-(m;\l)$ and belongs to $\polya^W$. 
Let  $D(m;q) \in \D(\frak a, \Sigma,m)$ denote 
the associated differential operator.
Since the $\Theta$-spherical functions solve (\ref{eq:hypereq}), 
we have on $A^+$
\begin{equation} \label{eq:eigen}
  D(m;q) \pedtheta{\varphi}(m;\l,a)=q(m;\l)\pedtheta{\varphi}(m;\l,a).
\end{equation}
The differential operator $D(m;q)$ might be singular on the set 
$\{\Delta=0\}$, but 
we can always choose $k \geq m$ so that $D_q:=\Delta^k D(m;q)$ is a 
$W$-invariant differential operator on $A$ with real analytic coefficients. 
Morover, (\ref{eq:eigen}) gives, for some constant $\kappa$
\begin{equation*}
  D_q\big(\pedtheta{\mathcal I}g\big)(a)=\kappa \Delta^{k-m}(a) 
D_m \mathcal F_A^{-1} (q g)(a) 
\end{equation*}
with $qg \in PW(C)$. The classical Paley-Wiener theorem now ensures that
$\supp  D_q (\pedtheta{\mathcal I}g) \subset \exp C$. 
We have already proven that $\pedtheta{\mathcal I}g$ is compactly supported.
If the leading symbol of $D_q$ never vanished, then Holmgren's uniqueness 
theorem would imply that 
\begin{equation}
\label{eq:holm}
\conv\big( \supp  D_q (\pedtheta{\mathcal I}g)\big)=
\conv\big( \supp  \pedtheta{\mathcal I}g\big).
\end{equation}
The problem is that, by construction, our differential operator
  $D_q$ has zeros along hyperplanes determined by $\a=0$ with 
$\a \in \Sigma$.
The necessary extension of Holmgren's uniqueness theorem was
accomplished in \cite{OP4}. Hence (\ref{eq:holm}) in fact holds.
For details we 
refer the reader to \cite{OP3} and \cite{OP4}. 

\appendix\section{$K_\e$-symmetric spaces with even multiplicities}
\label{app:classification}

In this appendix we report the infinitesimal classification of 
$K_\e$-symmetric spaces with even multiplicities by listing
the $K_\e$-symmetric pairs $(\frak g,\frak h)$ 
with even multiplicities for which $\frak g$ is simple
and noncompact. The list has been extracted 
from the classification due to Oshima and Sekiguchi \cite{OS}.
It is presented in three tables respectively collecting 
(for the even multiplicity case) the Riemannian symmetric pairs (Table 1),
the non-compactly causal (NCC) symmetric pairs (Table 2) and the other 
$K_\e$ symmetric pairs (Table 3). A non-Riemannian $K_\e$-symmetric pair
is said to be of type $K_\e I$ if its signature $\e$ comes from a gradation
of first kind according to \cite{Kane}. Otherwise it is said to be 
of type $K_\e II$. The symmetric pairs of type $K_\e I$ coincide with the 
NCC symmetric pairs. Table 3 therefore collects all symmetric pairs 
with even multiplicities of type $K_\e II$. 

The restricted root system $\Sigma$ of a $K_\e$-symmetric pair 
with even multiplicities has at most two root lengths. The classification
below shows that all multiplicities $m_\a$ of $\Sigma$ are equal, 
and moreover that they are all equal to $2$ for symmetric pairs of type
$K_\e II$.   
The restricted root system and multiplicities of a $K_\e$-symmetric pair 
$(\frak g, \frak h)$ coincide with
those of the corresponding Riemannian dual symmetric pair 
$(\frak g,\frak k)$. They are explicitly reported in Tables 
2 and 3 for the reader's convenience.

If $\Sigma$ is of type $X_n$ (with $X_n \in \{A_n,B_n,C_n,\dots\}$), then 
the index $n$ denotes the real rank of $\frak g$. The range for $n$ 
is chosen to avoid overlappings due to isomorphisms of symmetric spaces.
These isomorphisms arise from isomorphisms of the lower dimensional 
complex Lie algebras. We refer to \cite{He1}, Ch. X, \S 6, for more 
information. After each table below we report the relevant
symmetric pair isomorphisms.

\begin{table}[h]
\label{table:Riem}
\setlength{\extrarowheight}{4pt}
{\begin{tabular}{|>{$}c<{$}|>{$}c<{$}|>{$}c<{$}|>{$}c<{$}|>{$}c<{$}|>{$}c<{$}|}
\hline
\frak{g} &\frak{h}=\frak{k}&\Sigma&m_\alpha
\rule[-0.2cm]{0cm}{0.0cm}
&\\
\hline
\hline
\frak{sl}(n,\mathbb{C})& \frak{su}(n)
&A_{n-1}&2&n\geq 2
\rule[-0.2cm]{0cm}{0.0cm}
\cr\hline
\frak{so}(2n+1,\mathbb{C})&\frak{so}(2n+1)
&B_n&2&n\geq 2
\rule[-0.2cm]{0cm}{0.0cm}
\cr\hline
\frak{sp}(n,\mathbb{C})&\frak{sp}(n)
&C_n&2&n\geq 3
\rule[-0.2cm]{0cm}{0.0cm}
\cr\hline
\frak{so}(2n,\mathbb{C})&\frak{so}(2n)
&D_n&2&n\geq 4
\rule[-0.2cm]{0cm}{0.0cm}
\cr\hline
(\frak{e}_6)_{\smC}&\frak{e}_6&E_6&2
\rule[-0.2cm]{0cm}{0.0cm}
&\cr\hline
(\frak{e}_7)_{\smC}&\frak{e}_7&E_7&2
\rule[-0.2cm]{0cm}{0.0cm}
&\cr\hline
(\frak{e}_8)_{\smC}&\frak{e}_8&E_8&2
\rule[-0.2cm]{0cm}{0.0cm}
&\cr\hline
(\frak{f}_4)_{\smC}&\frak{f}_4&F_4&2
\rule[-0.2cm]{0cm}{0.0cm}
&\cr\hline
(\frak{g}_2)_{\smC}&\frak{g}_2&G_2&2
\rule[-0.2cm]{0cm}{0.0cm}
&\cr\hline
\frak{su}^*(2n)&\frak{sp}(n)&A_{n-1}&4&n\geq 2
\rule[-0.2cm]{0cm}{0.0cm}
\cr\hline
\frak{e}_{6(-26)}&\frak{f}_{4(-20)}&A_2&8
\rule[-0.2cm]{0cm}{0.0cm}
&\cr\hline
\frak{so}(2n+1,1)&\frak{so}(2n+1)&A_1&2n&n\geq 3
\rule[-0.2cm]{0cm}{0.0cm}
\cr\hline
\end{tabular}
\bigskip
\caption{Riemannian symmetric pairs with even multiplicities.} }
\end{table}

Special isomorphisms of Riemannian symmetric spaces with even multiplicities
are
\begin{alignat*}{2}
 \frak{so}(3,\C)=\frak{sp}(1,\C) &\approx \frak{sl}(2,\C), 
  &\quad
 \frak{so}(3)=\frak{sp}(1) &\approx \frak{su}(2)\,;\\
 \frak{sp}(2,\C) &\approx \frak{so}(5,\C),
  &\quad
\frak{sp}(2) &\approx \frak{so}(5)\,;\\
 \frak{so}(6,\C) &\approx \frak{sl}(4,\C),
  &\quad
\frak{so}(6) &\approx \frak{su}(4)\,;\\
 \frak{so}(3,1) &\approx \frak{sl}(2,\C),
  &\quad
\frak{so}(3) &\approx \frak{su}(2)\,;\\
 \frak{so}(5,1) &\approx \frak{su}^*(4),
  &\quad
\frak{so}(5) &\approx \frak{sp}(2)\,.
\end{alignat*}

The Lie algebra $\frak{so}(2,\C)$ is not semisimple.  
Observe  also that 
$\frak{so}(4,\C)\cong \frak{sl}(2,\C) \times \frak{sl}(2,\C)$ is not simple.
The structure of its homogeneous spaces can be therefore deduced from 
the structure of the homogeneous spaces of $\frak{sl}(2,\C)$.
 
\bigskip 
In the next table we list all the non-compactly causal, or $K_\e I$, 
symmetric pairs with even multiplicities. The third column reports the
subalgebra of $\frak g$ fixed by $\theta\theta_\e$, where $\theta_e$ is the 
involution associated with the $K_\e$-pair $(\frak g,\frak h)$.

\begin{table}[h]\label{table:NCC}
\setlength{\extrarowheight}{4pt}
\begin{tabular}{|>{$}c<{$}|>{$}c<{$}|>{$}c<{$}|>{$}c<{$}|>{$}c<{$}|>{$}c<{$}|}
\hline
\frak{g} &\frak{h}& \frak{g}^{\theta\theta_\e}&\Sigma& m_\alpha 
\rule[-0.2cm]{0cm}{0.0cm}
& 
 \\
\hline
\hline
\frak{sl}(n,\mathbb{C})&\frak{su}(n-j,j)
&\frak{sl}(n-j,\mathbb{C})\times\frak{sl}(j,\mathbb{C})\times\mathbb{C}
&A_{n-1}&2
&n\geq 2,\, 1\leq j\leq [n/2]
\rule[-0.2cm]{0cm}{0.0cm}
\cr\hline
\frak{so}(2n+1,\mathbb{C})&\frak{so}(2n-1,2)
&\frak{so}(2n-1,\mathbb{C})\times \mathbb{C}&B_n&2
&n\geq 2
\rule[-0.2cm]{0cm}{0.0cm}
\cr\hline
\frak{sp}(n,\mathbb{C})&\frak{sp}(n,\mathbb{R})
&\frak{gl}(n,\mathbb{C})&C_n&2
&n\geq 3
\rule[-0.2cm]{0cm}{0.0cm}
\cr\hline
\frak{so}(2n,\mathbb{C})&\frak{so}(2n-2,2)
&\frak{so}(2n-2,\mathbb{C})\times \mathbb{C}&D_n&2
&n\geq 4
\rule[-0.2cm]{0cm}{0.0cm}
\cr\hline
\frak{so}(2n,\mathbb{C})&\frak{so}^*(2n)
&\frak{gl}(n,\mathbb{C})&D_n&2
&n\geq 5
\rule[-0.2cm]{0cm}{0.0cm}
\cr\hline
(\frak{e}_6)_{\smC}&\frak{e}_{6(-14)}
&\frak{so}(10,\mathbb{C})\times \mathbb{C}&E_6&2
\rule[-0.2cm]{0cm}{0.0cm}
& 
\cr\hline
(\frak{e}_7)_{\smC}&\frak{e}_{7(-25)}
&(\frak{e}_6)_{\smC}\times\mathbb{C}&E_7&2
\rule[-0.2cm]{0cm}{0.0cm}
& 
\cr\hline
\frak{su}^*(2n)&\frak{sp}(n-j,j)
&\frak{su}^*(2(n-j))\times \frak{su}^*(2j)\times\mathbb{R}&A_{n-1}&4
&n\geq 2,\, 1\leq j \leq [n/2]
\rule[-0.2cm]{0cm}{0.0cm}
\cr\hline
\frak{e}_{6(-26)}&\frak{f}_{4(-20)}
&\frak{so}(9,1)\times \mathbb{R}&A_2&8&\cr\hline
\frak{so}(2n+1,1)&\frak{so}(2n,1)
&\frak{so}(2n+1)\times \mathbb{R}&A_1&2n
&n\geq 3
\rule[-0.2cm]{0cm}{0.0cm}
\cr\hline
\end{tabular}
\bigskip
\caption{Non-compactly causal symmetric pairs with even multiplicities.}  
\end{table}

\addtolength{\textheight}{40pt}

Special isomorphisms of NCC symmetric pairs with even multiplicities are
\begin{alignat*}{2}
 \frak{so}(3,\C)=\frak{sp}(1,\C) &\approx \frak{sl}(2,\C), 
  &\quad
 \frak{so}(1,2)\approx\frak{sp}(1,\R) &\approx \frak{su}(1,1)\,;\\
 \frak{sp}(2,\C) &\approx \frak{so}(5,\C),
  &\quad
\frak{sp}(2,\R) &\approx \frak{so}(3,2)\,;\\
 \frak{so}(6,\C) &\approx \frak{sl}(4,\C),
  &\quad
\frak{so}(4,2) &\approx \frak{su}(2,2)\,;\\
\frak{so}(6,\C) &\approx \frak{sl}(4,\C),
  &\quad
\frak{so}^*(6) &\approx \frak{su}(3,1)\,;\\
 \frak{so}(3,1) &\approx \frak{sl}(2,\C),
  &\quad
\frak{so}(2,1) &\approx \frak{su}(1,1)\,;\\
 \frak{so}(5,1) &\approx \frak{su}^*(4),
  &\quad
\frak{so}(4,1) &\approx \frak{sp}(1,1)\,;\\
 \frak{so}(8,\C) &= \frak{so} (8,\C),
  &\quad
\frak{so^*}(8) &\approx \frak{so}(2,6)\,;
\end{alignat*}
 
\smallskip 
The last table contains all 
the other $K_\e$-symmetric pairs, i.e. those of
type $K_\e II$, with even multiplicities.

\begin{table}[h] \label{table:nonNCC}
\setlength{\extrarowheight}{4pt}
{\begin{tabular}{|>{$}c<{$}|>{$}c<{$}|>{$}c<{$}|>{$}c<{$}|>{$}c<{$}|>{$}c<{$}|}
\hline
\frak{g} &\frak{h}& \frak{g}^{\theta\theta_\e}&\Sigma& m_\alpha 
\rule[-0.2cm]{0cm}{0.0cm}
& \\
\hline
\hline
\frak{so}(2n+1,\mathbb{C})&\frak{so}(2(n-j)+1,2j)
&\frak{so}(2(n-j)+1,\mathbb{C})\times
\frak{so}(2j,\mathbb{C})&B_n&2
&n \geq 2, \, 2\leq j\leq n
\rule[-0.2cm]{0cm}{0.0cm}
\cr\hline
\frak{sp}(n,\mathbb{C})&\frak{sp}(n-j,j)
&\frak{sp}(n-j,\mathbb{C})\times \frak{sp}(j,\mathbb{C})&C_n&2
&n \geq 3,\, 1\leq j \leq [n/2]
\rule[-0.2cm]{0cm}{0.0cm}
\cr\hline
\frak{so}(2n,\mathbb{C})&\frak{so}(2(n-j),2j)
&\frak{so}(2(n-j),\mathbb{C})\times\frak{so}(2j,\mathbb{C})&D_n&2
&n\geq 4,\, 2 \leq j \leq [n/2]
\rule[-0.2cm]{0cm}{0.0cm}
\cr\hline
(\frak{e}_6)_{\mathbb{C}}&\frak{e}_{6(2)}
&\frak{sl}(6,\mathbb{C})\times\frak{sl}(2,\mathbb{C})&E_6&2&
\cr\hline
(\frak{e}_7)_{\mathbb{C}}&\frak{e}_{7(7)}
&\frak{sl}(8,\mathbb{C})&E_7&2
\rule[-0.2cm]{0cm}{0.0cm}
&
\cr\hline
(\frak{e}_7)_{\mathbb{C}}&\frak{e}_{7(-5)}
&\frak{so}(12,\mathbb{C})\times\frak{sl}(2,\mathbb{C})&E_7&2
\rule[-0.2cm]{0cm}{0.0cm}
&
\cr\hline
(\frak{e}_8)_{\mathbb{C}}&\frak{e}_{8(8)}
&\frak{so}(16,\mathbb{C})&E_8&2
\rule[-0.2cm]{0cm}{0.0cm}
&
\cr\hline
(\frak{e}_8)_{\mathbb{C}}&\frak{e}_{8(-24)}
&e_{7\mathbb{C}}\times\frak{sl}(2,\mathbb{C})&E_8&2
\rule[-0.2cm]{0cm}{0.0cm}
&
\cr\hline
(\frak{f}_4)_{\mathbb{C}}&\frak{f}_{4(4)}
&\frak{sp}(3,\mathbb{C})\times\frak{sl}(2,\mathbb{C})&F_4&2
\rule[-0.2cm]{0cm}{0.0cm}
&
\cr\hline
(\frak{f}_4)_{\mathbb{C}}&\frak{f}_{4(-20)}
&\frak{so}(9,\mathbb{C}) &F_4&2&
\cr\hline
(\frak{g}_2)_{\smC}&\frak{g}_{2(2)}
&\frak{sl}(2,\mathbb{C})\times\frak{sl}(2,\mathbb{C})&G_2&2
\rule[-0.2cm]{0cm}{0.0cm}
&
\cr\hline
\end{tabular}
\bigskip
\caption{Other $K_\e$-symmetric pairs with even multiplicities.} }
\end{table}

\addtolength{\textheight}{-40pt}

A special isomorphism of $K_\e II$ symmetric pairs with even multiplicities
is
\begin{alignat*}{2}
 \frak{sp}(2,\C) &\approx \frak{so}(5,\C),
  &\quad
\frak{sp}(1,1) &\approx \frak{so}(1,4)\,.
\end{alignat*}


\begin{thebibliography}{mmBS00}

\bibitem[AO01]{AO01}
N.B. Andersen, and
G. {\'O}lafsson,
\newblock {\em A Paley-Wiener theorem
for the spherical Laplace transform
on causal symmetric spaces of
rank one},
\newblock  Proc. Amer. Math. Soc.
\textbf{129} (2001), 173--179


\bibitem[A{\'O}S00]{AOS}
N.B. Andersen, G.~{\'O}lafsson, and H.~Schlichtkrull,
\newblock {\em On the inversion of the Laplace and Abel transforms on causal
  symmetric spaces},
\newblock Preprint, 2000. 
\newblock To appear in {\em Forum Math.}

\bibitem[AU02]{AU}
N.B.~Andersen and J.M.~Unterberger,
\newblock {\em An application of shift operators to ordered
symmetric spaces}, 
Ann. Inst. Fourier (Grenoble) \textbf{52}  (2002),  no. 1, 275--288.


\bibitem[Bee88]{BeeComp}
R.~J.~Beerends,
\newblock {\em The Abel transform and shift operators},
\newblock Compositio Math. 66 (1988), no. 2, 145--197.

\bibitem[Bee87]{BeeCan}
\bysame,
\newblock {\em An introduction to the Abel transform},
\newblock  
Miniconference on harmonic analysis and operator algebras (Canberra, 1987), 
21--33, 
\newblock Proc. Centre Math. Anal. Austral. Nat. Univ., 15, Austral. 
Nat. Univ., Canberra, 1987.

\bibitem[vdBS93]{vdBSconv}
E.P. van den Ban and H.~Schlichtkrull,
\newblock {\em Convexity for invariant differential operators on semisimple
symmetric spaces}, 
\newblock {Compositio Math.}
   \textbf{89}  (1993), no. 3, 301--313.

\bibitem[vdBS97]{BS}
\bysame,
\newblock {\em The most-continuous part of the Plancherel decomposition
for a reductive symmetric space},
\newblock Ann. of Math. (2) \textbf{145} (1997), no. 2, 267--364.

\bibitem[Bou02]{Bou}
N. Bourbaki,
\newblock {\em Lie groups and Lie algebras},
\newblock Chapters 4--6,
\newblock Translated from the 1968 French original by Andrew Pressley.
\newblock   Elements of Mathematics. Springer-Verlag, Berlin, 2002.

\bibitem[Che91]{CherInv}
\newblock I.~Cherednik,
\newblock 
{\em A unification of Knizhnik-Zamolodchikov and Dunkl operators
via affine Hecke algebras},
\newblock  {Invent. Math.} \textbf{106} (1991), no. 2, 411--431.


\bibitem[Er53]{Er}
A.~Erd\'elyi et~al.
\newblock {\em {Higher transcendental functions}}, volume~1.
\newblock McGraw-Hill, New York, 1953.

\bibitem[Far95]{F95} J. Faraut, 
\newblock {\em Functions sph\'eriques sur un espace
sym\'etrique ordonn\'e de type Cayley}, Contemp Math.
\textbf{191} (1995), 41--55.

\bibitem[FH{\'O}94]{FHO}
J.~Faraut, J.~Hilgert, and G.~{\'O}lafsson,
\newblock {\em Spherical functions on ordered symmetric spaces},
\newblock {Ann. Inst. Fourier} \textbf{44} (1994), 927--966.      

\bibitem[Gan71]{G1}
R.~Gangolli,
\newblock {\em On the {P}lancherel formula and the {P}aley-{W}iener theorem for
  spherical functions on semisimple {L}ie groups},
\newblock {Ann. of Math. (2)} \textbf{93} (1971), 150--165.

\bibitem[Gan72]{G2}
\bysame,
\newblock Spherical functions on semisimple {L}ie groups.
\newblock In {\em Symmetric spaces (Short Courses, Washington Univ., St. Louis,
  Mo., 1969--1970)}, pages 41--92. Pure and Appl. Math., Vol. 8. Dekker, New
  York, 1972.

\bibitem[GV88]{GV}
R.~Gangolli and V.~S. Varadarajan,
\newblock {\em Harmonic analysis of spherical functions on real reductive
  groups}, 
\newblock  Ergebnisse der Mathematik
   und ihrer Grenzgebiete, 101.
\newblock Springer-Verlag, Berlin, 1988.

\bibitem[Gra97]{G97} P. Graczyk,
\newblock {\em Function c
on an ordered symmetric space},
\newblock Bull. Sci. Math. \textbf{121} (1997), 561--572.


\bibitem[Hba87]{Hba}
A.~Hba, 
\newblock {\em Analyse harmonique sur ${\rm SL}(3, H)$}.
\newblock C. R. Acad. Sci. Paris S\'er. I Math. 305 (1987), no. 3, 77--80.

\bibitem[Hec87]{Heck1}
G.~J. Heckman,
\newblock Root systems and hypergeometric functions. {I}{I}.
\newblock {\em Compositio Math.}, \textbf{64} (1987), no. 3, 353--373.

\bibitem[Heck91]{Heck2}
\bysame,
\newblock {\em An elementary approach to the hypergeometric
shift operators of Opdam},
\newblock Invent. Math.  \textbf{103}  (1991),  no. 2, 341--350.



\bibitem[Heck97]{HeckBou}
\bysame,
\newblock {\em Dunkl operators},
\newblock {Ast\'erisque} \textbf{245}, Exp.\ no.\ 828, 4, 223--246, 1997,
\newblock S\'eminaire Bourbaki, Vol.\ 1996/97.

\bibitem[HO87]{HOpd1}
G.~J. Heckman and E.~M. Opdam,
\newblock {\em Root systems and hypergeometric functions}. {I},
\newblock {Compositio Math.} \textbf{64} (1987), no. 3, 329--352.       

\bibitem[HS94]{HS}
G.~J. Heckman and H.~Schlichtkrull,
\newblock {\em {Harmonic analysis and special functions on symmetric spaces}},
\newblock Perspectives in Mathematics, 16.
\newblock Academic Press, Inc., San Diego, CA, 1994.

\bibitem[Hel66]{HePW}
S.~Helgason,
\newblock {\em An analogue of the Paley-Wiener theorem for the 
Fourier transform on certain symmetric spaces},
\newblock Math. Ann. \textbf{165} (1966), 297--308.

\bibitem[Hel78]{He1}
\bysame,
\newblock {\em {Differential geometry, Lie groups, and symmetric spaces}},
\newblock  Pure and Applied Mathematics, 80.
\newblock Academic Press, Inc., New York-London, 1978.

\bibitem[Hel84]{He2}
\bysame,
\newblock {\em {Groups and geometric analysis. Integral geometry, invariant
  differential operators, and spherical functions}},
\newblock Pure and Applied Mathematics, 113. 
\newblock Academic Press, Inc., Orlando, FL, 1984.

\bibitem[Hel94]{He3}
\bysame,
\newblock {\em {Geometric analysis on symmetric spaces}},
\newblock Mathematical Surveys and Monographs, 39. 
\newblock American Mathematical  Society, Providence, RI, 1994.



\bibitem[H{\'O}96]{HObook}
J.~Hilgert and G.~{\'O}lafsson,
\newblock {\em {Causal symmetric spaces. Geometry and Harmonic Analysis}},
\newblock  Perspectives in Mathematics, 18.
\newblock Academic Press, Inc., San Diego, CA, 1997.

\bibitem[Hoer90]{HoerLPDO}
L. H\"ormander,
\newblock {\em {The analysis of linear partial differential operators. I.
Distribution theory and Fourier analysis}},
\newblock  Springer Study Edition. Springer-Verlag, Berlin, 1990.

\bibitem[Hu90]{Hu}
J.E.~Humphreys, 
\newblock {\em {Reflection groups and Coxeter groups}},
\newblock  Cambridge Studies in Advanced Mathematics, 29. 
\newblock Cambridge University Press, Cambridge, 1990.

\bibitem[dJ93]{deJeu}
\newblock M. F. E. de Jeu, 
\newblock {\em  The Dunkl transform},
\newblock Invent. Math. \textbf{113}  (1993),  no. 1, 147--162.


\bibitem[JL01]{JL}
J. Jorgenson and S. Lang,
\newblock {\em {Spherical inversion on SL${\sb n}(R)$}},
\newblock Springer Monographs in Mathematics.
\newblock Springer-Verlag, New York, 2001.

\bibitem[Ka96]{Kane} 
S.~Kaneyuki,
\newblock {\em Signatures of roots and a new characterization of 
causal symmetric spaces}, 
\newblock In {\em {Topics in geometry}}, 213--229, 
Progr. Nonlinear Differential Equations Appl., 20, 
Birkhäuser, 1996. 

\bibitem[Koo84]{Koo84}
T.~H.~Koornwinder, 
\newblock {\em Jacobi functions and analysis on noncompact semisimple
Lie groups}. 
\newblock In: {\em Special functions: group theoretical aspects and
applications},  1--85, Math. Appl., Reidel, Dordrecht, 1984.

\bibitem[K{\'O}02]{KO02}
B.~Kr{\"o}tz and G.~{\'O}lafsson,
\newblock {\em The $c$-function for non-compactly causal symmetric spaces},
\newblock  {Invent. Math.} \textbf{149} (2002),  no. 3, 647--659.

\bibitem[K{\'O}03]{KO03}
\bysame,
\newblock {\em The $c$-function
for non-compactly causal symmetric
spaces and its relations to harmonic
analysis and representation
theory}, Preprint 2003. Submitted


\bibitem[{\'O}l97]{O97}
G.~{\'O}lafsson,
\newblock {Spherical functions and spherical Laplace transform on ordered
  symmetric spaces},
\newblock Preprint. Available at~{\tt
  http:\!//www.math.lsu.edu/$\null_{\widetilde{\null}}$preprint}, 1997.

\bibitem[{\'O}P01]{OP1}
G.~{\'O}lafsson and A.~Pasquale,
\newblock {\em On the meromorphic extension of the spherical functions on
  noncompactly causal symmetric spaces},
\newblock J. Funct. Anal. \textbf{181} (2001), no. 2, 346--401.

\bibitem[{\'O}P02]{OP2}
\bysame,
\newblock {\em Regularity properties of generalized Harish-Chandra expansions}.
\newblock In A. Strasburger et al. (eds.), 
{\em {Geometry and analysis on finite- and infinite-dimensional Lie groups}},
\newblock Banach Center Publications \textbf{55} (2002), 335--348.

\bibitem[{\'O}P03]{OP3}
\bysame,
\newblock {\em A Paley-Wiener theorem for the $\Theta$-spherical transform:
The even multiplicity case}. Preprint 2003, math.FA/0304361

\bibitem[{\'O}P04]{OP4}
\bysame,
\newblock {\em Support properties and Holmgren's uniqueness theorem
for invariant differential operators with hyperplane singularities.}
\newblock In preparation.

\bibitem[Opd88a]{Opd1}
E.M.~Opdam,
\newblock {\em Root systems and hypergeometric functions. III},
\newblock {Compositio Math.} \textbf{67} (1988), no. 1, 21--49. 

\bibitem[Opd88b]{Opd2}
\bysame,
\newblock {\em Root systems and hypergeometric functions. IV},
\newblock {Compositio Math.} \textbf{67} (1988), no. 2, 191--209.

\bibitem[Opd93]{OpdGauss}
\bysame,
\newblock {\em An analogue of the Gauss summation formula for
hypergeometric functions related to root systems},
\newblock {Math. Z.} \textbf{212} (1993), 313--336.

\bibitem[Opd95]{OpdActa}
\bysame,
\newblock {\em Harmonic analysis for certain representations of graded Hecke
  algebras}.
\newblock {Acta math.} \textbf{175} (1995), 75--121.

\bibitem[Opd00]{Opd00}
\bysame,
\newblock {\em Lecture notes on {D}unkl operators for real and complex
  reflection groups},
\newblock MSJ Memoirs, 8.
\newblock Mathematical Society of Japan, Tokyo, 2000.     

\bibitem[Opd89]{OpdShift}
\bysame,
\newblock {\em Some applications of hypergeometric shift operators},
\newblock {Invent. Math.} \textbf{98} (1989), 1--18.

\bibitem[OS80]{OS}
T.~Oshima and J.~Sekiguchi,
\newblock {\em Eigenspaces of invariant differential operators on an
affine symmetric space},
\newblock {Invent. Math.} \textbf{57} (1980), no. 1, 1--81.

\bibitem[Pa01]{PHab}
A.~Pasquale,
\newblock {\em A theory of $\Theta$-spherical functions},
\newblock Habilitationsschrift, Technische Universit\"at Clausthal, 2002.

\bibitem[Pa02a]{P1}
\bysame,
\newblock {\em A theory of $\Theta$-spherical functions},
\newblock Preprint, 2002. Submitted.

\bibitem[Pa02b]{P2}
\bysame,
\newblock {\em The $\Theta$-spherical transform and its inversion},
\newblock Preprint, 2002. To appear on {\em Math. Scand.}

\bibitem[Ros77]{Ros}
J.~Rosenberg,
\newblock {\em A quick proof of Harish-Chandra's Plancherel 
theorem for spherical functions on a semisimple Lie group},
\newblock Proc. Amer. Math. Soc. \textbf{63}  (1977), no. 1, 143--149.



\bibitem[Vre84]{Vre}
L.~Vretare,
\newblock {\em Formulas for elementary spherical functions and 
generalized Jacobi polynomials}, 
\newblock  {SIAM J. Math. Anal. } \textbf{15} (1984), 805--833.

\end{thebibliography}
\end{document}